\documentclass[12pt,a4paper]{article}
\usepackage{amscd,amsmath,amssymb,amsfonts,times,ascmac}

\usepackage{mathrsfs}
\usepackage[dvips]{color} 
\usepackage[all]{xy}
\usepackage{enumerate}
 \usepackage{fancybox} 
 \usepackage{tikz-cd}

\numberwithin{equation}{section}
    \newtheorem{thm}{Theorem}[section]
    \newtheorem{lem}[thm]{Lemma}
    
    \newtheorem{prop}[thm]{Proposition}
    \newtheorem{cor}[thm]{Corollary}

    \newtheorem{conj}[thm]{Conjecture}

    \newtheorem{rem}[thm]{Remark}
\DeclareMathAlphabet{\mathpzc}{OT1}{pzc}{m}{it}
\newcommand{\qed}
{\mbox{}\nolinebreak$\square$\medbreak\par}
\newenvironment{pf}{\par\smallskip\noindent\emph{Proof.}}{\hfill\qed\par\smallskip}
\newenvironment{pf*}[1]{\par\smallskip\noindent\emph{#1.}}{\hfill\qed\par\smallskip}
\newcommand{\bysame}{\hskip.3em \leavevmode\rule[.5ex]{2.5em}{.3pt}\hskip0.5em}
\begin{document}
\title{A generalization of the Ross symbols in higher $K$-groups and hypergeometric functions II}\author{M. Asakura
\footnote{
Department of Mathematics, Faculty of Sciences,
Hokkaido University, 
Sapporo 060-0810, JAPAN. 
\texttt{asakura@math.sci.hokudai.ac.jp}
}}
\date\empty
\maketitle


\def\can{\mathrm{can}}
\def\Gal{\mathrm{Gal}}
\def\ch{{\mathrm{ch}}}
\def\Coker{\mathrm{Coker}}
\def\crys{\mathrm{crys}}
\def\dlog{{\mathrm{dlog}}}
\def\dR{{\mathrm{d\hspace{-0.2pt}R}}}            
\def\et{{\mathrm{\acute{e}t}}}  
\def\Frac{{\mathrm{Frac}}}
\def\phami{\phantom{-}}
\def\id{{\mathrm{id}}}              
\def\Image{{\mathrm{Im}}}        
\def\Hom{{\mathrm{Hom}}}  
\def\Ext{{\mathrm{Ext}}}
\def\MHS{{\mathrm{MHS}}}  
  
\def\can{\mathrm{can}}
\def\arg{\mathrm{arg}}
\def\ch{{\mathrm{ch}}}
\def\Coker{\mathrm{Coker}}
\def\crys{\mathrm{crys}}
\def\dlog{d{\mathrm{log}}}
\def\dR{{\mathrm{d\hspace{-0.2pt}R}}}            
\def\et{{\mathrm{\acute{e}t}}}  
\def\Frac{\operatorname{Frac}}
\def\phami{\phantom{-}}
\def\id{{\mathrm{id}}}              
\def\Image{{\mathrm{Im}}}        
\def\Hom{\operatorname{Hom}}  
\def\Ext{{\mathrm{Ext}}}
\def\MHS{{\mathrm{MHS}}}  
  
\def\ker{\operatorname{Ker}}          
\def\zar{\mathrm{zar}}
\def\Pic{{\mathrm{Pic}}}
\def\CH{{\mathrm{CH}}}
\def\NS{{\mathrm{NS}}}
\def\NF{{\mathrm{NF}}}
\def\an{{\text{\it an}}}
\def\End{\operatorname{End}}
\def\pr{{\mathrm{pr}}}
\def\red{{\mathrm{red}}}
\def\Proj{\operatorname{Proj}}
\def\ord{\operatorname{ord}}
\def\rig{{\mathrm{rig}}}
\def\reg{{\mathrm{reg}}}          %
\def\res{{\mathrm{res}}}          %
\def\Res{\operatorname{Res}}
\def\Spec{\operatorname{Spec}}     
\def\prim{{\mathrm{prim}}}
\def\syn{{\mathrm{syn}}}
\def\cont{{\mathrm{cont}}}
\def\ln{{\operatorname{ln}}}

\def\bA{{\mathbb A}}
\def\bC{{\mathbb C}}
\def\C{{\mathbb C}}
\def\G{{\mathbb G}}
\def\bE{{\mathbb E}}
\def\bF{{\mathbb F}}
\def\F{{\mathbb F}}
\def\bG{{\mathbb G}}
\def\bH{{\mathbb H}}
\def\bJ{{\mathbb J}}
\def\bL{{\mathbb L}}
\def\cL{{\mathscr L}}
\def\bN{{\mathbb N}}
\def\bP{{\mathbb P}}
\def\P{{\mathbb P}}
\def\bQ{{\mathbb Q}}
\def\Q{{\mathbb Q}}
\def\bR{{\mathbb R}}
\def\R{{\mathbb R}}
\def\bZ{{\mathbb Z}}
\def\Z{{\mathbb Z}}
\def\cH{{\mathscr H}}
\def\cD{{\mathscr D}}
\def\cE{{\mathscr E}}
\def\cF{{\mathscr F}}
\def\cU{{\mathscr U}}
\def\cX{{\mathscr X}}
\def\cY{{\mathscr Y}}
\def\O{{\mathscr O}}
\def\cR{{\mathscr R}}
\def\cS{{\mathscr S}}
\def\cV{{\mathscr V}}
\def\cX{{\mathscr X}}
\def\cM{{\mathscr M}}
%
\def\ve{\varepsilon}
\def\vG{\varGamma}
\def\vg{\varGamma}
%
%
%
%
\def\lra{\longrightarrow}
\def\lla{\longleftarrow}
\def\Lra{\Longrightarrow}
\def\hra{\hookrightarrow}
\def\lmt{\longmapsto}
\def\ot{\otimes}
\def\op{\oplus}
\def\l{\lambda}
\def\Isoc{{\mathrm{Isoc}}}
\def\Fil{{\mathrm{Fil}}}
\def\Dw{{\mathrm{Dw}}}

\def\MHS{{\mathrm{MHS}}}
\def\MHM{{\mathrm{MHM}}}
\def\VMHS{{\mathrm{VMHS}}}
\def\sm{{\mathrm{sm}}}
\def\tr{{\mathrm{tr}}}
\def\fib{{\mathrm{fib}}}
\def\FIsoc{{F\text{-Isoc}}}
\def\FMIC{{F\text{-MIC}}}
\def\Log{{\mathscr{L}{og}}}
\def\FilFMIC{{\mathrm{Fil}\text{-}F\text{-}\mathrm{MIC}}}

\def\wt#1{\widetilde{#1}}
\def\wh#1{\widehat{#1}}
\def\spt{\sptilde}
\def\ol#1{\overline{#1}}
\def\ul#1{\underline{#1}}
\def\us#1#2{\underset{#1}{#2}}
\def\os#1#2{\overset{#1}{#2}}

\def\hF{{F_{\check{\underline{a}}}}}
\def\Ross{{\xi_{\mathrm{Ross}}}}
\def\HG{{\mathrm{HG}}}
\def\Gr{{\mathrm{Gr}}}
\def\bad{{\mathrm{bad}}}
\def\nur{{\omega_{i_0\ldots i_d}}}
\def\whnur{\widehat\omega_{i_0\ldots i_d}}
\def\ur{\eta_{i_0\ldots i_d}}
\def\whur{\widehat\eta_{i_0\ldots i_d}}

\def\Isoc{{\mathrm{Isoc}}}
\def\FIsoc{{F\text{-}\mathrm{Isoc}^\dag}}
\def\MIC{{\mathrm{MIC}}}
\def\FMIC{{F\text{-}\mathrm{MIC}}}
\def\FilMIC{{\mathrm{Fil}\text{-}\mathrm{MIC}}}
\def\rigsyn{{\mathrm{rig}\text{-}\mathrm{syn}}}
\def\FilFMIC{{\mathrm{Fil}\text{-}F\text{-}\mathrm{MIC}}}
\def\Log{{\mathscr{L}{og}}}
\def\Defn{\text{\Ovalbox{Def.}\quad}}
\newcommand{\sPol}{\mathop{\mathscr{P}\!\mathit{ol}}\nolimits}

\begin{abstract}
This is a sequel of the paper \cite{As-Ross1} where we introduced higher Ross symbols
in higher $K$-groups of the hypergeometric schemes, and discussed the Beilinson 
regulators.
In this paper we give its $p$-adic counterpart and an application to
the $p$-adic Beilinson conjecture for K3 surfaces of Picard number $20$.
\end{abstract}

\section{Introduction}
In \cite{As-Ross1}, we introduced a {\it higher Ross symbol}
\[
\Ross:=\left\{\frac{1-x_0}{1-\nu_0 x_0},\ldots,\frac{1-x_d}{1-\nu_d x_d}\right\}
\]
in the Milnor $K$-group of the affine ring of 
a {\it hypergeometric scheme}
\begin{equation}\label{intro-eq1}
U_t:(1-x_0^{n_0})\cdots(1-x_d^{n_d})=t
\end{equation} 
where $\nu_k$ is a $n_k$-th root of unity
(cf. \S \ref{rev-sect} and \S \ref{Rev-Ross-sect}).
This is a generalization of the Ross symbol $\{1-z,1-w\}$ in $K_2$ of the 
Fermat curve $z^n+w^n=1$ introduced in \cite{ross1}, \cite{ross2} (see \cite[4.1]{As-Ross1}
for the connection with our higher Ross symbols).
In \cite{As-Ross1} we discuss
the Beilinson regulator map (cf. \cite{schneider})
\[
\reg_B:K_{d+1}(U_t)^{(d+1)}\lra H^{d+1}_\cD(U_t,\Q(d+1))
\]
from Quillen's higher $K$-group to the Deligne-Beilinson cohomology
where $K_i(-)^{(j)}\subset K_i(-)\ot\Q$ denotes the Adams weight piece.
The main result \cite[Theorem 5.5]{As-Ross1} tells that
$\reg(\Ross)$ is a linear combination of complex analytic functions
\begin{equation}\label{intro-eq2}
\cF_{\ul a}(t):=\sum_{k=0}^d(\psi(a_k)+\gamma)+\log(t)
+a_0\cdots a_d\,t\, F_{\ul a}(t)
\end{equation}
where $\ul a=(a_0,\ldots,a_d)$ and 
\[
F_{\ul a}(t):={}_{d+3}F_{d+2}\left({a_0+1,\ldots,a_d+1,1,1\atop 2,\ldots,2};t\right),
\]
is the hypergeometric function (see \cite[5.1]{As-Ross1} for $\cF_{\ul a}(t)$, 
and \cite{slater} or \cite[15,16]{NIST} for
the general theory of hypergeometric functions).

\medskip

The purpose of this paper is to provide the $p$-adic counterpart of
\cite{As-Ross1},
namely we prove similar theorems in $p$-adic situation by replacing the subjects as follows,
\begin{align*}
\text{Beilinson regulator}&\rightsquigarrow\text{syntomic regulator}\\
\cF_{\ul a}(t)&\rightsquigarrow\cF^{(\sigma)}_{\ul a}(t)
\end{align*}
where $\cF^{(\sigma)}_{\ul a}(t)$ is a certain $p$-adic convergent function
introduced in \cite{New}. 
Let us explain it more precisely.
Let $W$ be the Witt ring of a perfect field of characteristic $p>0$, and
$K=\Frac W$ the fractional field.
For a smooth projective variety $X$ over $W$, we denote by
$H^\bullet_\syn(X,\Z_p(j))$ the syntomic cohomology of Fontaine-Messing
(cf. \cite[Chapter I]{Ka}).
More generally, let $U$ be a smooth $W$-scheme such that there is an embedding $U\hra X$ into a 
smooth projective $W$-scheme $X$ with $Z=X\setminus U$ 
a simple relative normal crossing divisor over $W$.
Then the log syntomic cohomology of $(X,Z)$ is defined (cf. \cite[\S 2.2]{Ts}),
which we denote by $H^\bullet_\syn(U,\Z_p(j))$.
In their recent paper \cite{NN}, Nekov\'a\v{r} and Niziol established 
the {\it syntomic regulator maps}
\[
\reg_\syn^{i,j}:K_i(U)\ot\Q\lra H^{2j-i}_\syn(U,\Q_p(j))
:=H^{2j-i}_\syn(U,\Z_p(j))\ot\Q
\]
from Quillen's algebraic $K$-groups.
These are the $p$-adic counterpart of the Beilinson regulator maps.
Let us take $U=U_\alpha$ the hypergeometric scheme \eqref{intro-eq1}
for $\alpha\in W$ and take the degrees $(i,j)=(d+1,d+1)$,
\begin{equation}\label{intro-reg}
\reg_\syn=\reg_\syn^{d+1,d+1}:K_{d+1}(U_\alpha)^{(d+1)}\lra
H^{d+1}_\syn(U_\alpha,\Q_p(d+1))\cong H^d_\dR(U_{\alpha,K}/K),
\end{equation}
where $U_{\alpha,K}:=U_\alpha\times_WK$.
We then discuss the element $\reg_\syn(\Ross)$ in the de Rham cohomology.
Although the author does not know whether
the higher Ross symbol $\Ross\in K_{d+1}(U_\alpha)$ lies in the image of
$K_{d+1}(X_\alpha)$ with $X_\alpha\supset U_\alpha$ a smooth compactification
(see \cite[4.2]{As-Ross1} for more details), one can show
that $\reg_\syn(\Ross)$ lies in $W_dH^d_\dR(U_{\alpha,K}/K)$
 (Lemma \ref{regsyn-lem1}),
where $W_\bullet H^*_\dR(U_{\alpha,K}/K)$ denotes the weight filtration by
Deligne.
The main theorem of this paper (=Theorem \ref{main-2}) describes
$\reg_\syn(\Ross)$ by a linear combination of 
\[\cF^{(\sigma)}_{\ul a}(t)\]
introduced in \cite[\S 2]{New},
which is the $p$-adic counterpart of $\cF_{\ul a}(t)$
(see \S \ref{log-sect} below for the review on $\cF^{(\sigma)}_{\ul a}(t)$).
The proof is based on the theory of $F$-isocrystals,
especially the main result in \cite{AM}, and also uses the congruence relation for
$\cF^{(\sigma)}_{\ul a}(t)$ proven in \cite[3.2]{New}.

Our main theorem has an application to the study of the $p$-adic Beilinson conjecture
by Perrin-Riou.
Conceptually saying, the conjecture
 asserts that the special values of the $p$-adic
$L$-functions are described by the syntomic regulators up to $\Q^\times$.
There are several previous works by many people (\cite{BK}, \cite{BD}, \cite{KLZ},
\cite{Niklas}, etc.). We also refer the recent paper
\cite{AC} where a number of numerical verifications for $K_2$ of elliptic curves over $\Q$
are given apart from the method of the Beilinson-Kato elements.
In this paper we shall discuss the $p$-adic Beilinson conjecture
for singular K3 surfaces over $\Q$ (``singular" means the Picard number $20$).
Let $h^2_\tr(X)=h^2(X)/\NS(X)$ denote the transcendental motive (cf. \cite[7.2.2]{KMP}).
If $X$ is a singular K3 surface over $\Q$, it is 2-dimansional, and
there is a Hecke eigenform $f$ of weight $3$ with complex multiplication
such that
$L(h^2_\tr(X),s)=L(f,s)$ by a theorem of Livn\'e \cite{Li}.
The $p$-adic Beilinson conjecture for $K_3(X)^{(3)}$ 
is formulated as the relation between
the syntomic regulator and the special value of the $p$-adic $L$-function $L_p(f,\chi,s)$.
See Conjecture \ref{pBconj} for the precise statement.
Our hypergeometric schemes provide several explicit examples 
of singular K3 surfaces.
The hypergeometric scheme
\[
(1-x_0^2)(1-x_1^2)(1-x_2^2)=\alpha,\quad\alpha\in\Q\setminus\{0,1\}
\]
is a K3 surface with the Picard number $\geq 19$.
This is isogenous to the K3 surface
\[
w^2=u_1u_2(1+u_1)(1+u_2)(u_1-\alpha u_2),
\] 
studied by Ahlgren, Ono and Penniston \cite{ono}.
A complete list of $\alpha$'s such that the Picard number 20 (i.e. singular K3) is known. 
For such an example, one can employ the higher Ross symbol together with
our main theorem (=Theorem \ref{main-2}), so that we have a formulation of
the $p$-adic Beilinson conjecture in terms of our $p$-adic hypergeometric functions
$\cF^{(\sigma)}_{\ul a}(t)$'s (Conjectures \ref{pRZconj}, \ref{pRZconj-1}).
The advantage of our formulation is that it allows the numerical verifications,
while we have not succeeded a theoretical proof.
In \S \ref{K3-1-sect}, we construct other examples of singular K3 surfaces arising from
the hypergeometric schemes, and give a description of the $p$-adic regulators
in terms of $\cF^{(\sigma)}_{\ul a}(t)$ (Theorem \ref{ellK3-thm}).

\section{$p$-adic Hypergeometric Functions}
\subsection{Hypergeometric series}
Let $K$ be a field of characteristic zero.
Let $(\alpha)_n=\alpha
(\alpha+1)\cdots(\alpha+n-1)$ be the Pochhammer symbol.
For $\ul a=(a_0,\ldots,a_d)\in K^{d+1}$, the power series
\begin{equation}\label{HG-eq1}
F_{\ul a}(t)={}_{d+1}F_d\left({a_0,\ldots,a_d\atop 1,\ldots,1};t\right)=
\sum_{n=0}^\infty\frac{(a_0)_n}{n!}\cdots\frac{(a_d)_n}{n!}t^n
\end{equation}
is called the {\it hypergeometric series}.
In case $d=1$, this is also referred to as the Gaussian hypergeometric series.
If $a_i\not\in \Z$ for all $i$,
this is characterized as the unique solution (up to scalar) 
of the {\it hypergeometric differential operator}
\begin{equation}\label{HG-eq2}
P_{\HG,\ul a}:=D^{d+1}-t(D+a_0)\cdots(D+a_d),\quad D:=t\frac{d}{dt}.
\end{equation}
in the ring $K[[t,t^{-1}]]$ of Laurent power series.
\subsection{Dwork's Hypergeometric Functions}\label{Dw-sect}
Let $p$ be a prime.
For $\ul a=(a_0,\ldots,a_d)\in \Z_p^{d+1}$, the hypergeometric series
\[
F_{\ul a}(t)
=\sum_{n=0}^\infty\frac{(a_0)_n}{n!}\cdots\frac{(a_d)_n}{n!}t^n
\]
has coefficients in $\Z_p$.
In his paper \cite{Dwork-p-cycle},
Dwork discovered that a certain ratio of hypergeometric series 
is a uniform limit of a sequence of rational functions.

For $\alpha\in \Z_p$, 
let $\alpha'$ denote the Dwork prime, which is defined to be $(\alpha+k)/p$
where $k\in \{0,1,\ldots,p-1\}$ such that $\alpha+k\equiv 0$ mod $p$.
Define the $i$-th Dwork prime by
$\alpha^{(i)}=(\alpha^{(i-1)})'$ and $\alpha^{(0)}:=\alpha$.
Write $\ul a'=(a'_1,\ldots,a'_s)$ and 
$\ul a^{(i)}=(a^{(i)}_1,\ldots,a^{(i)}_s)$.
{\it Dwork's $p$-adic hypergeometric function} is defined to be a power series
\[
\cF^\Dw_{\ul a}(t):=F_{\ul a}(t)/F_{\ul a'}(t^p)\in \Z_p[[t]].
\]
Let $W=W(\ol\F_p)$ be the Witt ring of $\ol\F_p$.
Let $c\in 1+pW$ and $\sigma$ a $p$-th Frobenius on $W[[t]]$ 
given by $\sigma(t)=ct^p$, $c\in 1+pW$.
A slight modification of $\cF^\Dw_{\ul a}(t)$ is
\begin{flalign*}
&\Defn\hspace{2cm} 
\cF^{\Dw,\sigma}_{\ul a}(t):=F_{\ul a}(t)/F_{\ul a'}(t^{\sigma})\in W[[t]].&
\end{flalign*}
Dwork discovered certain congruence relations 
which we refer to as the {\it Dwork congruence}.
The precise statement is as follows.
For a power series $f(t)=\sum_{i\geq 0}a_it^i$,
we denote
$[f(t)]_{<k}=\sum_{0\leq i<k}a_it^i$ the truncated polynomial.
\begin{thm}[{\cite[p.37, Thm. 2, p.45]{Dwork-p-cycle}}]\label{D-C-eq0}
For any $n\geq 1$, we have
\[
\cF^{\Dw,\sigma}_{\ul a}(t)\equiv \frac{[F_{\ul a}(t)]_{<p^n}}{[F_{\ul a'}(t^{\sigma})]_{<p^n}}
\mod p^nW[[t]].
\]
\end{thm}
The Dwork congruence implies that $\cF^{\Dw,\sigma}_{\ul a}(t)$
is a $p$-adic analytic function in the following way.
Put
\begin{equation}\label{D-C-h-eq}
h_{\ul a}(t):=\prod_{i=0}^N[F_{\ul a^{(i)}}(t)]_{<p}
\end{equation}
with $N$ sufficiently large such that
$
\{\ol{[F_{\ul a^{(i)}}(t)]}_{<p}\}_{i\geq0}
=\{\ol{[F_{\ul a^{(i)}}(t)]}_{<p}\}_{0\leq i\leq N}$ as subsets of $\F_p[t]$
where $\ol{f(t)}:=f(t)$ mod $p$. 
Let
\begin{equation}\label{D-C-eq1}
\Z_p[t,h_{\ul a}(t)^{-1}]^\wedge:=\varprojlim_{n\geq1}
\bigg(\Z_p/p^n\Z_p[t,h_{\ul a}(t)^{-1}]\bigg)
\end{equation}
be the $p$-adic completion. The ring does not depend on the choice of $N$.
An immediate consequence of the Dwork congruence is
\[
\cF^{\Dw,\sigma}_{\ul a}(t)\in W\ot_{\Z_p}\Z_p[t,h_{\ul a}(t)^{-1}]^\wedge.
\]
As an another application of 
the Dwork congruence, one can show
\begin{equation}\label{D-C-eq2}
\frac{\frac{d^j}{dt^j}F_{\ul a}(t)}{F_{\ul a}(t)}\equiv
\frac{\frac{d^j}{dt^j}[F_{\ul a}(t)]_{<p^n}}{[F_{\ul a}(t)]_{<p^n}}
\mod p^n\Z_p[[t]]
\end{equation}
for all $n\geq 1$ and $j\geq0$ in the same way as the proof of 
\cite[(3.14)]{Dwork-p-cycle}.

\subsection{$p$-adic Hypergeometric Functions of logarithmic type \cite{New}}
\label{log-sect}
In \cite{New}, we introduce a certain new $p$-adic hypergeometric function, 
which is different from Dwork's one.
We recall it here.
Define a continuous function
\begin{equation}\label{wt-polygamma-def}
\wt{\psi}_p(z):=\lim_{n\in\Z_{>0},n\to z}\sum_{1\leq k<n,\,p\nmid k}
\frac{1}{k}
\end{equation}
on $\Z_p$
where ``$n\to z$" means the limit with respect to the $p$-adic metric (\cite[2.2]{New}).
Let $c\in 1+p W$ and $\sigma(t)=ct^p$ be as before.
Let
\[
G_{\ul a}^{(\sigma)}(t):=\sum_{i=0}^d\wt\psi_p(a_i)
-p^{-1}\log(c)+
\int_0^t(F_{\underline{a}}(t)-F_{\ul a'}(t^{\sigma}))\frac{dt}{t}
\]
where $\log(z)$ is the Iwasawa logarithmic function.
Then we define
\begin{flalign*}
&\Defn\hspace{2cm} 
\cF^{(\sigma)}_{\underline{a}}(t):=G^{(\sigma)}_{\ul a}(t)/F_{\underline{a}}(t),&
\end{flalign*}
and call the {\it $p$-adic hypergeometric functions of logarithmic type}.
This is a power series with coefficients in $W$.

\medskip

There are congruence relations that are similar to the Dwork congruence
(Theorem \ref{D-C-eq0}).
\begin{thm}[{\cite[Theorem 3.3]{New}}]\label{log-thm1}
If $c\in 1+2pW$, then there are congruence relations
\[
\cF^{(\sigma)}_{\ul a}(t)\equiv \frac{[G^{(\sigma)}_{\ul a}(t)]_{<p^n}}{[F_{\ul a'}(t^\sigma)]_{<p^n}}
\mod p^nW[[t]],\quad n\geq 1.
\]
If $p=2$ and $c\in 1+2W$, then the congruence holds modulo $p^{n-1}W[[t]]$.
Hence $\cF^{(\sigma)}_{\ul a}(t)$ belongs to the $p$-adic completion 
$W[t,h_{\ul a}(t)^{-1}]^\wedge$ of the ring $W[t,h_{\ul a}(t)^{-1}]$.
\end{thm}
The above congruence plays a key role in the proof of the main theorem (Theorem \ref{main-1}).
\section{Hypergeometric Schemes}
\subsection{Review from \cite{As-Ross1}}\label{rev-sect}
Let $A$ be a commutative ring.
Let $d\geq 1$ and $n_i\geq 2$ $(i=0,1,\ldots,d)$ be integers such that
$n_0\cdots n_d$ is invertible in $A$.
We call an affine scheme
\begin{equation}\label{HG-eq0}
U=\Spec A[x_0,\ldots,x_d]/((1-x_0^{n_0})\cdots(1-x_d^{n_d})-t),\quad t\in A
\end{equation}
the {\it hypergeometric scheme} over $A$ (\cite[\S 2.1]{As-Ross1}).
It is easy to see that $U$ is smooth over $A$ if $t(1-t)\in A^\times$.
\begin{prop}\label{compact-prop}
Assume that 
$t(1-t)\in A^\times$ and $A$ is an integral domain.
Then 
there is an open immersion $U\hra X$ into a smooth projective $A$-scheme $X$ 
such that $Z=X\setminus U$ is a 
relative simple NCD over $A$. 
\end{prop}
\begin{pf}
\cite[Proposition 2.1]{As-Ross1}.
\end{pf}
Hereafter we assume that $A$ is an integral domain and
$t(1-t)\in A^\times$.
Let $\mu_n=\mu_n(A)$ denote the group of $n$-th roots of unity in $A$.
A finite abelian group $G=\mu_{n_0}\times \cdots\times \mu_{n_d}$ acts on 
$U$ in a way that
$(x_0,\ldots,x_d)\mapsto (\nu_0x_0,\ldots,\nu_dx_d)$
for $\ul \nu=(\nu_0,\ldots,\nu_d)\in G$.
In this way the de Rham cohomology groups $H^\bullet_\dR(U/A)$ are endowed with
the structure of $A[G]$-modules.
For a $A[G]$-module $H$ and $(i_0,\ldots,i_d)\in \Z^{d+1}$
we denote by
\begin{equation}\label{rev-eigen-eq}
H(i_0,\ldots,i_d)=\{x\in H\mid\ul \nu x=\nu_0^{i_0}\cdots\nu_d^{i_d}x,\,\forall\ul\nu\in G\}. 
\end{equation}
the simultaneous eigenspace.

Let $A=K[t,(t-t^2)^{-1}]$ with $K$ a field of characteristic zero.
Suppose that $K$ contains primitive $n_i$-th roots of unity for all $i$.
We denote by $W_\bullet H^i_\dR(U/A)$ the weight filtration. 
In particular $W_i H^i_\dR(U/A)=\Image[H^i_\dR(X/A)\to H^i_\dR(U/A)]$.
Put $I:=\{(i_0,\ldots,i_d)\in\Z^{d+1}\mid 0\leq i_k<n_k\}$
and $I_+:=\{(i_0,\ldots,i_d)\in\Z^{d+1}\mid 0< i_k<n_k\}$.
There is the decomposition
\[
H^\bullet_\dR(U/A)=\bigoplus_{(i_0,\ldots,i_d)\in I}H^\bullet_\dR(U/A)(i_0,\ldots,i_d)
\]
of the de Rham cohomology group.
\begin{thm}[{\cite[\S 3.2 Summary (1), (2)]{As-Ross1}}]\label{dim-cor1}
\begin{enumerate}
\item[$(1)$]
The $A$-submodule
\[
\bigoplus_{(i_0,\ldots,i_d)\in I\setminus I_+}H^\bullet_\dR(U/A)(i_0,\ldots,i_d)
\]
is generated by 
exterior products of
$\{dx_k/(x_k-\nu_k)\mid \nu_k\in \mu_{n_k}\}$.
\item[$(2)$]
If $i<d$, then
\[
\bigoplus_{(i_0,\ldots,i_d)\in I_+}H^i_\dR(U/A)(i_0,\ldots,i_d)=0.
\]
In particular $W_iH^i_\dR(U/A)=0$ for $0<i<d$ by this and {\rm (1)}.
\item[$(3)$]
\[
\bigoplus_{(i_0,\ldots,i_d)\in I_+}H^d_\dR(U/A)(i_0,\ldots,i_d)=W_dH^d_\dR(U/A).
\]
\end{enumerate}
\end{thm}

For each $(i_0,\ldots,i_d)\in I_+$, we put, cf. \cite[\S 2.4]{As-Ross1} 
\begin{flalign}\label{form-eq1}
&\Defn\hspace{2cm}
\omega_{i_0\ldots i_d}:=n_0^{-1}x_0^{i_0-n_0}x_1^{i_1-1}\cdots x_d^{i_d-1}
\frac{dx_1\cdots dx_d}{(1-x_1^{n_1})\cdots(1-x_d^{n_d})}&
\end{flalign}

a regular $d$-form in $\vg(X,\Omega^d_{X/A})(i_0,\ldots,i_d)$.
\begin{thm}[{\cite[Corollary 3.6]{As-Ross1}}]\label{dim-cor2}
Let $\cD:=K\langle t,(t-t^2)^{-1},\frac{d}{dt}\rangle$ be the Weyl algebra of $\Spec A$, which acts on $H^d_\dR(U/A)$.
Let $(i_0,\ldots,i_d)\in I_+$ and put $a_k:=1-i_k/n_k$.
Let $P_{\HG,\ul a}$ be the hypergeometric differential operator \eqref{HG-eq2}. 
Then $P_{\HG,\ul a}(\omega_{i_0\ldots i_d})=0$ and the homomorphism
\[
\cD/\cD P_{\HG,\ul a}\os{\cong}{\lra}
H^d_\dR(U/A)(i_0,\ldots,i_d),\quad P\longmapsto P(\omega_{i_0\ldots i_d})
\]
of $\cD$-modules is bijective. In particular, this is an irreducible $\cD$-module.
\end{thm}
In view of Theorem \ref{dim-cor2}, we think
the piece $H^d_\dR(U/A)(i_0,\ldots,i_d)$ of being
an {\it hypergeometric motive} associated to the hypergeometric function \eqref{HG-eq1}.

\bigskip

\begin{itembox}[l]{\bf Convention}
For an integer $n>0$, let $\Z_{(n)}=S^{-1}\Z$ 
denote the ring of fractions with respect to the multiplicative set 
$S=\{s\in\Z\mid \gcd(s,n)=1\}$. 
For $j\in \Z_{(n)}$, let $[j]_n$ denote the unique integer such that
$0\leq [j]_n<n$ and $[j]_n\equiv j$ 
mod $n\Z_{(n)}$.
We then extend the notation as follows,
\[
H(i_0,\ldots,i_d):=H([i_0]_{n_0},\ldots,[i_d]_{n_d}),\quad
(i_0,\ldots, i_d)\in \prod_{k=0}^d \Z_{(n_k)},
\]
\[
\omega_{i_0\ldots i_d}:=\omega_{[i_0]_{n_0}\ldots[i_d]_{n_d}},\quad
(i_0,\ldots, i_d)\in \prod_{k=0}^d (\Z_{(n_k)}\setminus n_k\Z_{(n_k)}).
\]
\end{itembox}
\subsection{Pairing $Q$ on $W_dH^d_\dR(U/A)$}\label{pairing-sect}
Let $K$ be a field of characteristic zero, and $A$ an integral smooth $K$-algebra.
Let $U$
be the hypergeometric scheme \eqref{HG-eq0} over $A$ of relative dimension $d$
(we do not assume that $K$ contains primitive $n_i$-th roots of unity).
We construct a natural pairing 
\begin{equation}\label{pairing-eq}
Q:W_dH^d_\dR(U/A)\ot_AW_dH^d_\dR(U/A)\lra A
\end{equation}
in the following way.
Let $X\supset U$ be a compactification such that $X\to\Spec A$ is smooth projective
with connected fibers (e.g. Proposition \ref{compact-prop}).
We fix an ample class $[\omega]\in H^2_\dR(X/A)$, and let $L$ be
the Lefschetz operator on $H^\bullet_\dR(X/A)$
i.e. the cup-product with $[\omega]$.
Let
\[
H^d_\dR(X/A)=\bigoplus_{i\geq 0}L^iH^{d-2i}_\dR(X/A)_\prim
\]
be the Lefschetz decomposition by the primitive parts of cohomology.
Let $j^*:H^\bullet_\dR(X/A)\to H^\bullet_\dR(U/A)$ be the pull-back by
the immersion $j:U\hra X$.
Then, for $i>0$, the image of the component $L^iH^{d-2i}_\dR(X/A)_\prim$
vanishes by Theorem \ref{dim-cor1} (2), so that the map
\[
j^*:H^d_\dR(X/A)_\prim\lra W_dH^d_\dR(U/A)
\]
is surjective. Let $T\subset H^d_\dR(X/A)_\prim$ be the kernel, and 
$T^\perp\subset H^d_\dR(X/A)_\prim$ the orthogonal complement with respect to
 the cup-product pairing $Q_X$ on $H^d_\dR(X/A)_\prim$. 
Since $Q_X$ induces the polarization on the primitive part $H^d_\dR(X/A)_\prim$, 
the restricted pairings on $T$ and $T^\perp$ are non-degenerate (e.g. \cite[Corollary 2.12]{MHS}),
and hence $T^\perp\cong W_dH^d_\dR(U/A)$.
Then we define the pairing \eqref{pairing-eq} to be the induced one from $Q_X$
restricted on $T^\perp$.

\begin{prop}\label{pairing-prop}
\begin{itemize}
\item[ \rm\textbf{(Q1)}]
The pairing $Q$ is $(-1)^d$-symmetric and non-degenerate,
\item[\rm \textbf{(Q2)}]
Suppose that $K$ contains primitive $n_i$-th roots of unity for all $i$. Then
$Q(\sigma x,\sigma y)=Q(x,y)$ for $\sigma\in G$.
\item[ \rm\textbf{(Q3)}]
$Q(\theta x,y)+Q(x,\theta y)=\theta(Q(x,y))$ for any
derivative $\theta$ on $A$.
\item[ \rm\textbf{(Q4)}]
$Q(F^p,F^q)=0$ if $p+q>d$ and
$Q(F^p,F^q)=A$ if $p+q=d$ and $p,q\geq0$, where $F^\bullet$ is the Hodge filtration.
\end{itemize}
When $K$ contains primitive $n_i$-th roots of unity for all $i$,
the property {\bf(Q2)} implies that $Q$ induces a perfect pairing
\begin{equation}\label{pairing-eq6}
Q:W_dH^d_\dR(U/A)(i_0,\ldots,i_d)\ot_A W_dH^d_\dR(U/A)(-i_0,\ldots,-i_d)\lra A.
\end{equation}
\end{prop}
 \begin{pf}
Everything but \textbf{(Q2)}
is obvious from the construction.
We should be careful of \textbf{(Q2)}
as we do not assume that the action of $G$ extends on $X$. In proving \textbf{(Q2)}, 
we give an alternative construction of $Q$.
 We may replace $K$ with $\ol K$.
Let $H^\bullet_{\dR,c}(U/A)$ denote the de Rham cohomology with compact support, on which $G$ acts.
There is the cup-product
\[
H^d_{\dR,c}(U/A)\ot_A H^d_{\dR}(U/A)\lra 
H^{2d}_{\dR,c}(U/A)\os{\cong}{\to}H^{2d}_\dR(X/A)\cong A,
\]
which induces a pairing
\[
Q_c:\Gr^W_dH^d_{\dR,c}(U/A)\ot_A \Gr^W_dH^d_{\dR,c}(U/A)\lra A.
\]
This is
compatible with $Q_X$ under the natural map $\Gr^W_dH^d_{\dR,c}(U/A)
\to H^d_\dR(X/A)$, and satisfies $Q_c(\sigma x,\sigma y)=Q_c(x,y)$ for $\sigma\in G$.
We note that $W_dH^d_{\dR,c}(U/A)=H^d_{\dR,c}(U/A)$.
If one can show that 
\[
\Image(H^d_{\dR,c}(U/A))\subset H^d_\dR(X/A)_\prim
\]
and that the composition 
\[
u:H^d_{\dR,c}(U/A)\to H^d_\dR(X/A)\to W_dH^d_\dR(U/A)
\]
is surjective, then the pairing induced from $Q_c$ agrees with $Q$, and hence \textbf{(Q2)} follows.
The former is equivalent to that
the composition
$H^d_{\dR,c}(U/A)\to H^d_\dR(X/A)
\os{L}{\to} H^{d+2}(X/A)$ is zero.
This agrees with the composition
$H^d_{\dR,c}(U/A)\os{L_U}{\to} H^{d+2}_{\dR,c}(U/A)\to H^{d+2}_\dR(X/A)$
where $L_U$ is the cup-product with $[\omega]|_U\in W_2H^2_\dR(U/A)$.
Therefore it is enough to show  $[\omega]|_U=0$.
If $d\ne2$, this follows from the vanishing $W_2H^2_\dR(X/A)=0$ (Theorem \ref{dim-cor1} (2)).
If $d=2$, it follows from Theorem \ref{dim-cor1} (3) and Theorem \ref{dim-cor2}
 that there is no non-zero element $z\in W_2H^2_\dR(U/A)$ such that $Dz=0$, and in particular one has
 $[\omega]|_U=0$.

We show that $u$ is surjective.
This is equivalent to 
 the surjectivity of the natural map
\[
u_{i_0\ldots i_d}:H^d_{\dR,c}(U/A)(i_0,\ldots,i_d)\lra W_dH^d_\dR(U/A)(i_0,\ldots,i_d)
\]
for each $(i_0,\ldots,i_d)\in I_+$ by Theorem \ref{dim-cor1} (3).
Since this is a homomorphism of $\cD$-modules, and the right 
hand side is irreducible by Theorem \ref{dim-cor2},
it is enough to see the non-vanishing $u_{i_0\ldots i_d}\ne0$.
To see this, it is enough to show that
\[
\xymatrix{
F^dH^d_{\dR,c}(U/A)\ar[r]&F^dH^d_{\dR}(X/A)\ar[r]^-\cong& F^dW_dH^d_\dR(U/A)
}
\]
is surjective.
However this follows from
the exact sequence
\[
H^d_{\dR,c}(U/A)\lra H^d_{\dR}(X/A)\lra H^d_{\dR}(Z/A)
\]
where $Z=X\setminus U$, and the vanishing $F^dH^d_\dR(Z/A)=0$ as $\dim(Z/A)=d-1$.
This completes the proof.
\end{pf}
\begin{lem}
The pairing $Q$ does not depend on the choice of $X$ and $[\omega]$.
\end{lem}
\begin{pf}
The alternative construction of $Q$ in the proof of Proposition \ref{pairing-prop} does not depend on
$X$ and $[\omega]$.
\end{pf}
\section{Unit root formula for Hypergeometric Schemes}
Throughout this section,
let $W=W(\F)$ be the Witt ring of an algebraically closed field $\F$
of characteristic $p>0$, $K:=\Frac(W)$ the fractional field
and put $A:=W[t,(t-t^2)^{-1}]$.

\subsection{Unit Root Vectors}\label{unit-sect}
Let $n_0,\ldots,n_d>1$ be integers and $p$ a prime such that
$p\nmid n_0\cdots n_d$.
Let 
\[
U=\Spec A[x_0,\ldots,x_d]/((1-x_0^{n_0})\cdots(1-x_d^{n_d})-t)
\]
be the hypergeometric scheme over $A$. 
We write $A_K:=K[t,(t-t^2)^{-1}]$, and
$U_K:=U\times_AA_K$.

Let $\cD:=K\langle t,(t-t^2)^{-1},\frac{d}{dt}\rangle$ be the Weyl algebra.
Put
$
 D:=t\frac{d}{dt}
$.
Let $(i_0,\ldots,i_d)\in \prod_{k=0}^d\Z_{(n_k)}$ satisfy
$[i_k]_{n_k}\ne0$ for all $k$ 
(see {\bf Convention} in \S \ref{rev-sect} for the notation).
Recall from Theorems \ref{dim-cor1} and \ref{dim-cor2} the eigenspace
\[
H^d_\dR(U_K/A_K)(i_0,\ldots,i_d)
=W_dH^d_\dR(U_K/A_K)(i_0,\ldots,i_d)
=\sum _{k=0}^d A_KD^k\omega_{i_0\ldots i_d},
\]
which is stable under the action of $\cD$.
We write
$H_{i_0\ldots i_d}(U_K/A_K)=H^d_\dR(U_K/A_K)(i_0,\ldots,i_d)$
for simplicity of notation.
\begin{prop}\label{hat-unit-prop}
Let $H_{i_0\ldots i_d}(U_K/A_K)_{K((t))}:=
K((t))\ot_{A_K} H_{i_0\ldots i_d}(U_K/A_K)$ on which the action of $\cD$ extends 
in a natural way.
Put $a_k:=1-[i_k]_{n_k}/n_k$. 
Let $s_k\in \Q$ be defined by
$(t+a_0)\cdots(t+a_d)=t^{d+1}+s_1t^d+\cdots+s_{d+1}$, and 
put $q_{d-m}:=-s_{m+1}t/(1-t)$
for $m=0,1,\ldots,d$. Put $\check{\ul a}:=(1-a_0,\ldots,1-a_d)$ and
\[
y_d:=(1-t)\hF(t)=(1-t){}_{d+1}F_d\left({1-a_0,\ldots,1-a_d\atop 1,\ldots,1};t\right).
\] 
For $0\leq i<d$
define $y_i$ by
$y_i+Dy_{i+1}=q_{i+1}y_d$.
Put
\begin{equation}\label{unit-eq1}
\whur:=y_0\nur+y_1D\nur+\cdots+y_dD^d\nur\in H_{i_0\ldots i_d}(U_K/A_K)_{K((t))}.
\end{equation}
Then
\begin{equation}\label{unit-eq2}
\ker[D:H_{i_0\ldots i_d}(U_K/A_K)_{K((t))}\to H_{i_0\ldots i_d}(U_K/A_K)_{K((t))}]
=K\whur.
\end{equation}
\end{prop}
\begin{pf}
Recall Theorem \ref{dim-cor2}. 
$H_{i_0\ldots i_d}(U_K/A_K)_{K((t))}$ is a free $K((t))$-module with
basis $\{D^k\nur\mid k=0,1,\ldots,d\}$ and the differential operator 
$P_{\HG,\ul a}=D^{d+1}-t(D+a_0)\cdots(D+a_d)=(1-t)(D^{d+1}+q_dD^d+\cdots+q_0)$
annihilates $\nur$.
Therefore
\begin{align*}
D\left(\sum_{k=0}^dz_kD^k\nur\right)&=\sum_{k=0}^dD(z_k)D^{k}\nur+z_kD^{k+1}\nur\\
&=\sum_{k=1}^d(z_{k-1}+D(z_k))D^{k}\nur+D(z_0)\nur+z_dD^{d+1}\nur\\
&=\sum_{k=1}^d(z_{k-1}+D(z_k)-q_kz_d)D^{k}\nur+(D(z_0)-q_dz_d)\nur
\end{align*}
vanishes if and only if
\begin{equation}\label{unit-lem1-eq1}
z_i+D(z_{i+1})=q_{i+1}z_d\,(0\leq i\leq d-1),\quad D(z_0)=q_dz_d.
\end{equation}
Put a differential operator
\[
P:=D^{d+1}-D^d\star q_d+\cdots+(-1)^dD \star q_1+(-1)^{d+1}q_0
\]
 where $\star$ denotes
the composition of operators to make distinctions between $D\star f\in \cD$ and 
$D(f)\in K((t))$.
Then \eqref{unit-lem1-eq1} is equivalent to
\[
z_i+D(z_{i+1})=q_{i+1}z_d\,(0\leq i\leq d-1),\quad P(z_d)=0,
\]
so that the assertion is reduced to show that $y_d=(1-t)\hF(t)$ is the
unique solution (up to scalar) in $K((t))$ of the differential equation $P(y)=0$.
One has
\begin{align*}
P\star (1-t)&=D^{d+1}\star (1-t)-\sum_{m=0}^{d}(-1)^m\left(\sum_{i_0<\cdots<i_m}a_{i_0}\cdots a_{i_m}\right)D^{d-m}\star t\\
&=D^{d+1}-t(D+1)^{d+1}-\sum_{m=0}^{d}(-1)^{m}\left(\sum_{i_0<\cdots<i_m}a_{i_0}\cdots a_{i_m}\right)t(D+1)^{d-m}\\
&=D^{d+1}-t(D+1-a_0)\cdots(D+1-a_d)\\
&=P_{\HG,\check{\ul a}}.
\end{align*}
Therefore $y_d=(1-t)\hF(t)$ is the unique solution for $P$.
\end{pf}

We define a {\it unit root vector}
\begin{flalign}\label{unit-eq3}
&\Defn
\ur:=\hF(t)^{-1}\whur=\frac{y_0}{\hF(t)}\nur+\frac{y_1}{\hF(t)}D\nur+\cdots+(1-t)D^d\nur.&
\end{flalign}

\begin{lem}\label{unit-lem1}
Let $h_{\ul a}(t)$ be the polynomial as in \eqref{D-C-h-eq}.
Then we have
\[
\frac{y_i}{\hF(t)}\in
\bigg(\Z_p[t,h_{\check{\ul a}}(t)^{-1}]^\wedge\bigg)[(1-t)^{-1}],\quad i=0,1,\ldots,d
\]
where $(-)^\wedge$ denotes the $p$-adic completion, and
hence 
\[
\ur\in 
A[h_{\check{\ul a}}(t)^{-1}]^\wedge\ot_AH_{i_0\ldots i_d}(U_K/A_K).
\]
\end{lem}
\begin{pf}
By the construction of $y_i$, they are linear combinations
of $D^j(y_d)$ over a ring $\Z_{(p)}[t,(1-t)^{-1}]$ and hence
one can write
\[
y_{i}
=\sum_j g_j\frac{d^j\hF(t)}{dt^j}
\]
by some $g_j\in \Z_{(p)}[t,(t-t^2)^{-1}]$.
Now
the assertion follows from \eqref{D-C-eq2}.
\end{pf}
Let $h(t):=\prod_ah_{\ul a}(t)$
where $\ul a=(i_0/n_0,\ldots,i_d/n_d)$ runs over all $(d+1)$-tuple of integers 
$(i_0,\ldots,i_d)$ such that $0<i_k<n_k$.
Put $B:=A[h(t)^{-1}]$, $B_K:=K\ot_WB$, and 
\begin{equation}\label{wh-B}
\wh B:=A[h(t)^{-1}]^\wedge=\varprojlim_{n}\bigg(W/p^nW[t,(t-t^2)^{-1},h(t)^{-1}]\bigg)
\end{equation}
the $p$-adic completion, and 
$\wh B_K:=K\ot_W \wh B$.
Thanks to Lemma \ref{unit-lem1}, the unit root vector
$\ur$ belongs to $H_{i_0\ldots i_d}(U_K/A_K)_{\wh B_K}:=\wh B_K\ot_{A_K}H_{i_0\ldots i_d}(U_K/A_K)$.
We call the direct summand
\begin{flalign*}
&\Defn\hspace{2cm}
H^{\text{unit}}_{i_0\ldots i_d}(U_K/A_K)_{\wh B_K}:=B_K\,\ur&
\end{flalign*}
of $H_{i_0\ldots i_d}(U_K/A_K)_{\wh B_K}$
the {\it unit root subspace}.
Recall from \S \ref{pairing-sect} the perfect pairing \eqref{pairing-eq6}
\begin{equation}\label{unit-eq6}
Q:H_{i_0\ldots i_d}(U_K/A_K)\ot_{A_K} H_{-i_0,\ldots, -i_d}(U_K/A_K)\lra A_K.
\end{equation}
Tensoring with $\wh B_K$, we have the pairing on the $\wh B_K$-modules, which
we also write by $Q$.
Define a $\wh B_K$-submodule
\begin{flalign}\label{N-defn}
&\Defn\hspace{2cm}
VH_{i_0\ldots i_d}(U_K/A_K)_{\wh B_K}\subset H_{i_0\ldots i_d}(U_K/A_K)_{\wh B_K}&
\qquad
\end{flalign}
to be the exact annihilator of the unit root part 
$H^{\text{unit}}_{-i_0,\ldots, -i_d}(U_K/A_K)_{\wh B_K}$.
By definition, the pairing 
\begin{equation}\label{unit-eq5}
H_{i_0\ldots i_d}(U_K/A_K)_{\wh B_K}/VH_{i_0\ldots i_d}(U_K/A_K)_{\wh B_K}
\ot
H^{\text{unit}}_{-i_0,\ldots, -i_d}(U_K/A_K)_{\wh B_K}
\lra \wh B_K
\end{equation}
is perfect. We shall later see $\ur\in VH_{i_0\ldots i_d}(U_K/A_K)_{\wh B_K}$ 
(Corollary \ref{RC-cor2}).
Let $W((t))^\wedge$ be the $p$-adic completion and $K((t))^\wedge:=K\ot_W
W((t))^\wedge$. 
We write $VH_{i_0\ldots i_d}(U_K/A_K)_{K((t))^\wedge}:=K((t))^\wedge\ot_{A_K}
VH_{i_0\ldots i_d}(U_K/A_K)$. Note that $\wh B_K\subset K((t))^\wedge$. 
\begin{lem}\label{unit-lem2}
\begin{enumerate}
\item[$(1)$]
$VH_{i_0\ldots i_d}(U_K/A_K)_{\wh B_K}$ is stable under the action of $\cD$.
\item[$(2)$]
$H_{i_0\ldots i_d}(U_K/A_K)_{\wh B_K}/VH_{i_0\ldots i_d}(U_K/A_K)_{\wh B_K}\cong
\wh B_K$ is generated by $\nur$.
\item[$(3)$]
$Q(\nur,\eta_{-i_0,\ldots, -i_d})
=C\in \Q^\times$.
\end{enumerate}
\end{lem}
\begin{pf}
(1) is immediate from the fact that $D\whur=0$ (Proposition \ref{hat-unit-prop}). 
To see (2), it is enough to
show $Q(\nur,\eta_{-i_0,\ldots, -i_d})\ne0$.
Write $\omega=\nur$ and $\check\omega=\omega_{-i_0,\ldots, -i_d}$.
Recall from \cite[Corollary 3.10]{As-Ross1} the fact that
$\Gr_F^iH_{i_0\ldots i_d}(U_K/A_K)$ is a free $A_K$-module
with basis $D^{d-i}\omega$. 
Therefore,
thanks to the property {\bf(Q4)} in Proposition \ref{pairing-prop}, we have
\begin{equation}\label{unit-lem2-eq3}
Q(D^i\omega,D^j\check\omega)
=0
\end{equation}
for any $(i,j)$ such that  $i+j<d$ and $i,j\geq0$.
In particular 
\begin{equation}\label{unit-lem2-eq3-1}
Q(\nur,\eta_{-i_0,\ldots, -i_d})=(1-t)
Q(\omega,D^d\check\omega).
\end{equation}
If $Q(\nur,\eta_{-i_0,\ldots, -i_d})=0$, then 
$Q(\omega,D^i\check\omega)$ vanishes for all $i\in\Z_{\geq0}$, which
contradicts with the fact that $Q$ is a perfect pairing.
Hence $Q(\nur,\eta_{-i_0,\ldots, -i_d})\ne0$.
We show (3).
Since
\[
0=DQ(D^i\omega,D^j\check\omega)=Q(D^{i+1}\omega,D^{j}\check\omega)
+Q(D^{i}\omega,D^{j+1}\check\omega)
\]
for $i+j=d-1$ and $i,j\geq0$ by \eqref{unit-lem2-eq3}, one has
\begin{equation}\label{unit-lem2-eq1}
(-1)^iQ(\omega,D^d\check\omega)=
Q(D^i\omega,D^{d-i}\check\omega),\quad i=0,1,\ldots,d.
\end{equation}
Applying $D$ on the both sides, one has
\[
(-1)^iDQ(\omega,D^d\check\omega)=
Q(D^{i+1}\omega,D^{d-i}\check\omega)
+Q(D^{i}\omega,D^{d-i+1}\check\omega).
\]
Taking the alternating sum of the both sides, one has
\begin{equation}\label{unit-lem2-eq2}
(d+1)DQ(\omega,D^d\check\omega)
=(-1)^dQ(D^{d+1}\omega,\check\omega)
+Q(\omega,D^{d+1}\check\omega).
\end{equation}
Using $P_{\HG,\ul a}\omega=0$ and $P_{\HG,\check{\ul a}}\check\omega=0$
together with \eqref{unit-lem2-eq3}, one has
\begin{align*}
Q(D^{d+1}\omega,\check\omega)
&=(a_0+\cdots+a_d)\frac{t}{1-t}Q(D^d\omega,\check\omega)\\
&=(-1)^d(a_0+\cdots+a_d)\frac{t}{1-t}Q(\omega,D^d\check\omega)\quad
\text{(by \eqref{unit-lem2-eq1})},\\
Q(\omega,D^{d+1}\check\omega)&=(d+1-(a_0+\cdots+a_d))
\frac{t}{1-t}Q(\omega,D^d\check\omega).
\end{align*}
Apply the above to \eqref{unit-lem2-eq2}, then
\[
(d+1)DQ(\omega,D^d\check\omega)
=(d+1)\frac{t}{1-t}Q(\omega,D^d\check\omega)
\quad\Longleftrightarrow\quad
\frac{d}{dt}
Q(\omega,D^d\check\omega)
=\frac{1}{1-t}Q(\omega,D^d\check\omega).
\]
This implies $Q(\omega,D^d\check\omega)=C(1-t)^{-1}$ with some $C\in K$.
Hence $Q(\nur,\eta_{-i_0,\ldots, -i_d})=C$ by \eqref{unit-lem2-eq3-1}, and this is not zero
by (2). 
Since the pairing $Q$ and the elements $\omega,\check\omega$ are defined over 
a ring $\Q(t)$, it turns out that
$C\in K^\times\cap \Q(t)=\Q^\times$.
This completes the proof of (3).
\end{pf}
\subsection{Rigid cohomology of Hypergeometric schemes}\label{RC-sect}
We write $X_K=X\times_WK$ and $X_\F=X\times_W\F$ for a $W$-scheme $X$.
Let $T$ be a smooth affine scheme over $W$.
Let $\O(T)^\dag$ denote the weak completion
of $\O(T)$, and write $\O(T)^\dag_K:=\O(T)^\dag\ot_WK$.
We fix a $p$-th Frobenius $\sigma$ on $\O(T)^\dag$, namely
it satisfies
$\sigma(x)\equiv x^p$ mod $p$, and
$\sigma$ on $W$ agrees with the $p$-th Frobenius on the Witt ring.
For a smooth
morphism $g:V\to T$, the $i$-th rigid cohomology group
\[
H^i_\rig(V_\F/T_\F)
:=\Gamma(T_K^\an, R^ig_{\rig}j^{\dag}_{V}\O_{V_K^\an}),
\]
is defined (cf. \cite[Definition 2.12]{AM}).
The $p$-th Frobenius $\Phi$ that is $\sigma$-linear is defined on
the rigid cohomology group in the canonical way. 
In this paper, we often employ the following fundamental fact.
\begin{thm}\label{coh-thm}
Suppose that either of the following conditions {\rm (1)} and {\rm (2)} holds.
\begin{itemize}
\item[\rm (1)]
The morphism $g:V\to T$ is smooth projective.
\item[\rm (2)]
There is a completion $T\hra \ol T$ into 
a smooth projective $W$-scheme $\ol T$ 
such that $E:=\ol T\setminus T$ is a smooth divisor over $W$,
and there are $\ol Y$, $Y$ and $V$ which fit into
ommutative squares
\[
\xymatrix{
V\ar[r]^j\ar[d]_g&Y\ar@{}[rd]|{\square}
\ar[d]_f\ar[r]^{j'}&\ol Y\ar[d]^{\bar f}\\
T\ar@{=}[r]&T\ar[r]^\subset&\ol T
}
\]
satisfying the following condistions, where $\square$ denotes a cartesian diagram.
\begin{itemize}
\item[\rm(2-1)]
$\ol Y$ is smooth projective over $W$.
The arrows $j$ and $j'$ are open embeddings, $\bar f$ is projective, and $f$ is
smooth projective. 
\item[\rm(2-2)]
The scheme-theoretic inverse image $D:=\bar f^{-1}(E)$ is a 
relative simple normal crossing divisor over $W$, and each multiplicity
is prime to $p$.
\item[\rm(2-3)]
$Z:=Y\setminus V$ is a relative simple normal crossing divisor over $T$.
\item[\rm(2-4)]
Let $\bar Z$ be the closure of $Z$ in $\bar Y$.
Then, $\bar Z+D$ is a 
relative simple normal crossing divisor over $W$.
\end{itemize}
\end{itemize}
Then the $i$-th relative rigid cohomology sheaf $R^ig_{\rig}j^{\dag}_V\O_{V_K^{\an}}$
is a coherent $j^{\dag}_T\O_{T_K^{\an}}$-module for each $i$.
Consequently, $H^i_\rig(V_\F/T_\F)$ is a locally free $\O(T)^\dag_K$-module
of finite rank, 
and the comparison map 
\[
c:\O(T)^\dag_K\ot_{\O(T_K)}H^i_\dR(V_K/T_K)\lra H^i_\rig(V_\F/T_\F)
\]
is bijective for each $i$
(see \cite[2.5]{AM} for the construction of the comparison map).
\end{thm}
\begin{pf}
See \cite[Propostions 2.15, 2.17]{AM} for the proof.
We note that the essential point is 
Kiehl's theorem for proper case and 
\cite[Theorem 2.2]{Shiho} for non-proper case.
\end{pf}
We turn to the setting in \S \ref{unit-sect}.
We write $A_\F:=\F[t,(t-t^2)^{-1}]$ and
$U_\F:=U\times_AA_\F$ etc.
For $c\in 1+pW$, let $\sigma$ be a $F$-linear $p$-th Frobenius 
on $A^\dag_K$ given by $\sigma(t)=c t^p$.
Let $X\supset U$ be the compactification as in Proposition \ref{compact-prop},
and $Z:=X\setminus U$ the complements which is a relative simple 
normal crossing divisor over $A$. 
We shall later see that $U\to\Spec  A$ satisfies the condition 
Theorem \ref{coh-thm} (2) if $p>\mathrm{lcm}(n_0,\ldots,n_d)$
(Lemma \ref{completion-lem}).
However in this section, we work under a weaker assumption
that $p\nmid n_0\cdots n_d$.

\medskip

Let $Z=\cup_k Z_k$ be the irreducible decomposition.
Then one has an exact sequence
\[
\bigoplus_k H^{i-2}_\dR(Z_{k,K}/A_K)\lra H^i_\dR(X_K/A_K)\lra
W_iH^i_\dR(U_K/A_K)\lra 0,
\]
so that one has the Frobenius structure on
\[
A_K^\dag\ot_{A_K}W_iH^i_\dR(U_K/A_K)
\]
compatible with the $\cD$-module structure thanks to  Theorem \ref{coh-thm}.
Write
\[
H_{i_0\ldots i_d}(U_K/A_K)_{A_K^\dag}:=
A_K^\dag\ot_{A_K}W_dH^d_\dR(U_K/A_K)(i_0,\ldots,i_d)
=A_K^\dag\ot_{A_K}H^d_\dR(U_K/A_K)(i_0,\ldots,i_d).
\]
for $(i_0,\ldots,i_d)\in \prod_{k=0}^d\Z_{(n_k)}$ satisfying $i_k\not\equiv0$ mod $n_k$.
Tensoring \eqref{unit-eq6} with $A_K^\dag$, we have a pairing
\[
H_{i_0\ldots i_d}(U_K/A_K)_{A^\dag_K}
\ot_{A^\dag_K} 
H_{-i_0\ldots -i_d}(U_K/A_K)_{A^\dag_K}
\lra A^\dag_K\ot_AH^{2d}_{c,\dR}(U_K/A_K)\cong A^\dag_K
\]
which we also write by $Q$.
Since the cup-product is compatible with the Frobenius, we have
\begin{description}
\item[(Q5)]
$Q(\Phi(x),\Phi(y))=p^d\sigma(Q(x,y))$.
\end{description}
The Frobenius $\sigma$ extends on $\wh B_K$ and $K((t))^\wedge$ in a natural way.
According to this, we extend $\Phi$ to that on $H_{i_0\ldots i_d}(U_K/A_K)_{\wh B_K}$
and $H_{i_0\ldots i_d}(U_K/A_K)_{K((t))^\wedge}
=K((t))^\wedge\ot_AH_{i_0\ldots i_d}(U_K/A_K)$, and
denote by the same notation.

\begin{lem}\label{RC-lem1}
Let $H_{i_0\ldots i_d}$ denote $H_{i_0\ldots i_d}(U_K/A_K)_R$ 
where $R$ is either of $A_K^\dag$, $\wh B_K$ or $K((t))^\wedge$, and
$VH_{i_0\ldots i_d}$ denotes $VH_{i_0\ldots i_d}(U_K/A_K)_R$ 
(see \eqref{N-defn} for the notation).
Then
\begin{enumerate}
\item[{\rm (i)}]
$\Phi(H_{p^{-1}i_0\ldots p^{-1}i_d})\subset H_{i_0\ldots i_d}$.
\item[{\rm (ii)}]
$\Phi(H_{p^{-1}i_0\ldots p^{-1}i_d}^{\mathrm{unit}})\subset H_{i_0\ldots i_d}^{\mathrm{unit}}$.
\item[{\rm (iii)}]
$\Phi(VH_{p^{-1}i_0\ldots p^{-1}i_d})\subset VH_{i_0\ldots i_d}$.
\end{enumerate}
\end{lem}
\begin{pf}
Since each $H_{i_0\ldots i_d}$ is an irreducible $\cD$-module, there is unique
$(j_0,\ldots,j_d)\in \prod_{k=0}^d\Z_{(n_k)}$ with $j_k\not\equiv0$ mod $n_k$
such that
\[
\Phi(H_{i_0\ldots i_d})\subset H_{j_0\ldots j_d}.
\]
The assertion (i) is equivalent to that $H_{j_0\ldots j_d}=H_{pi_0,\ldots, pi_d}$.
To show this,
we restrict the cohomology at a fiber $t=a$
where $a\in W$ such that $a\not\equiv 0,1$ mod $p$.
We take $ca^p=F(a)$ so that one has $\sigma(t)|_{t=a}
=F(a)$. Let $U_a:=U\times_A\Spec A/(t-a)$ and
$X_a:=X\times_A\Spec A/(t-a)$.
Then $\Phi$ induces the Frobenius
\[
\Phi_a(H^d_\dR(U_{a,K}/K)(i_0,\ldots,i_d))\subset H^d_\dR(U_{a,K}/K)(j_0,\ldots,j_d).
\]
which agrees with the Frobenius on $H^d_\rig(U_{a,\F}/\F)$ under
the isomorphism $H^\bullet_\dR(U_{a,K}/K)\cong H^\bullet_\rig(U_{a,\F}/\F)$
(Theorem \ref{coh-thm}).
Now the assertion
follows from the fact that $g\Phi_a=\Phi_a g$ on $H^\bullet_\rig(U_{a,\F}/\F)$
for $g\in G$
(see \S \ref{rev-sect} for the group $G$).
This completes the proof of (i).

We show (ii) and (iii).
Since $\ker D$ in \eqref{unit-eq2} is stable under the action of $\Phi$, we have
\begin{equation}\label{RC-lem1-eq1}
\Phi(\wh\eta_{p^{-1}i_0\ldots p^{-1}i_d})=C\whur
\end{equation}
for some $C\in K^\times$ ($C$ is not zero by {\bf (Q5)} on noticing that
$Q$ is a perfect pairing).
This implies
\begin{equation}\label{RC-lem1-eq2}
\Phi(\eta_{p^{-1}i_0\ldots p^{-1}i_d})
=C
\frac{F_{\check{\ul a}}(t)}{F_{\check{\ul a}^\prime}(t^{\sigma})}
\eta_{i_0\ldots i_d}
=C
\cF^{\Dw,\sigma}_{\check{\ul a}}(t)
\eta_{i_0\ldots i_d}
\end{equation}
where $\check{\ul a}=(1-a_0,\ldots,1-a_d)$ and 
$\check{\ul a}^\prime=(1-a^\prime_0,\ldots,1-a^\prime_d)$ denotes the Dwork prime.
Hence (ii) follows, and
(iii) is immediate from (ii) and the definition of $VH_{i_0\ldots i_d}$.
\end{pf}

\begin{thm}\label{RC-thm1}
Let $(i_0,\ldots,i_d)$ satisfy $i_k\not\equiv 0$ mod $n_k$ for all $k$.
Then 
\[
\Phi(\omega_{p^{-1}i_0\ldots p^{-1}i_d})\equiv p^d\cF^{\Dw,\sigma}_{\ul a}(t)^{-1}\nur
\mod VH_{i_0\ldots i_d}(U_K/A_K)_{\wh B_K}
\]
\end{thm}
\begin{pf}
Put $\wh\omega_{i_0\ldots i_d}:=F_{\ul a}(t)^{-1}\omega_{i_0\ldots i_d}$.
By Lemma \ref{unit-lem2} (3), 
$Q(\wh\omega_{i_0\ldots i_d},\wh\eta_{-i_0,\ldots,-i_d})$ is a non-zero constant.
By \eqref{RC-lem1-eq1} together with {\bf(Q5)}, we have
\[
\Phi(\wh\omega_{p^{-1}i_0\ldots p^{-1}i_d})\equiv
p^dC'\wh\omega_{i_0\ldots i_d}
\,\Leftrightarrow\,
\Phi(\omega_{p^{-1}i_0\ldots p^{-1}i_d})\equiv p^dC'\cF^{\Dw,\sigma}_{\ul a}(t)^{-1}\nur
\mod VH_{i_0\ldots i_d}
\]
with some $C'\in K^\times$.
We show $C'=1$. To do this, we employ the log crystalline cohomology.
We recall from \cite[Lemma 2.2]{As-Ross1} a cartesian diagram
\[
\xymatrix{
X\ar[d]\ar[r]&Y\ar[d]^f\\
\Spec A\ar[r]^-\subset&\Spec W[t]
}
\]
where $Y$ is smooth over $W$, and the scheme theoretic fiber $D=f^{-1}(O)$
with $O:=\Spec W[t]/(t)$ is reduced and relative simple NCD over $W$.
Let $\cY:=X\times_AW[[t]]\to \Spec W[[t]]$. 
We endow the log-structures on $\cY$ and $\Spec W[[t]]$ by the divisors $D$ and $O$
respectively, and let
\[
H^i_{\text{log-crys}}((\cY,D)/(W[[t]],O))
\]
be the $i$-th log crystalline cohomology group endowed with the $p$-th Frobenius $\Phi_{(\cY,D)}$ (\cite{Ka-log}).
The cohomology is described by the relative log de Rham complex.
Let
\[
\omega^\bullet_\cY:=\Omega^\bullet_{\cY/W}(\log D), \quad
\omega^\bullet_{W[[t]]}:=\Omega^\bullet_{W[[t]]/W}(\log O)
\] 
denote the log de Rham complex, and
\[
\omega^\bullet_{\cY/W[[t]]}:=\Coker\left[\frac{dt}{t}\ot\omega^{\bullet-1}_\cY\to \omega^\bullet_\cY\right].
\]
Then there is the canonical isomorphism
$H^i_{\text{log-crys}}((\cY,D)/(W[[t]],O))
\cong H^i_\zar(\cY,\omega_{\cY/W[[t]]}^\bullet)$
(\cite[Theorem 6.4]{Ka-log}).
We have the canonical homomorphism
\[
H^i_\rig(X_\F/A_\F)\cong A_K^\dag\ot_{A_K}H^i_\dR(X_K/A_K)
\lra
\Q\ot_\Z W((t))^\wedge\ot_{W[[t]]}H^i_\zar(\cY,\omega_{\cY/W[[t]]}^\bullet)
\]
where $(-)^\wedge$ denote the $p$-adic completion, and 
$\Phi$ and $\Phi_{\cY,D)}$ are compatible.
For $\ul\nu=(\nu_0,\ldots,\nu_d)\in G$, 
let $P_{\ul\nu}$ denote the subscheme of $D$ defined by
$\{x_0-\nu_0=\cdots=x_d-\nu_d=0\}$. 
Let $R$ be the composition
\[
\omega_{\cY/W[[t]]}^\bullet\os{\wedge\frac{dt}{t}}{\lra}\Omega^{\bullet+1}_{\cY/W}(\log D)\os{\Res}{\lra}
\bigoplus_{\ul\nu\in G} \O_{P_{\ul\nu}}[-d]
\]
of the complexes where $\Res$ is the Poincare residue map.
This induces
\[
\xymatrix{
H^d_{\text{log-crys}}((\cY,D)/(W[[t]],O))
\ar[r]^-R& \bigoplus_{\ul\nu\in G} W\cdot P_{\ul\nu}\\
H^d_\zar(\cY,\omega_{\cY/W[[t]]}^\bullet)\ar@{=}[u]
}\]
and it satisfies $R\circ\Phi_{(\cY,D)}=p^dF\circ R$.
We claim that the submodules $tH^d_\zar(\cY,\omega_{\cY/W[[t]]}^\bullet)$
and $DH^d_\zar(\cY,\omega_{\cY/W[[t]]}^\bullet)$ is annihilated by $R$ where $D:=t\frac{d}{dt}$.
The former is clear from the definition of the Poincare residue. To see the latter,
we recall that $D$ is defined to be the connecting homomorphism arising
from 
\[
0\lra\frac{dt}{t}\ot\omega^{\bullet-1}_{\cY/W[[t]]}\lra\omega^\bullet_{\cY}\lra\omega^\bullet_{\cY/W[[t]]}\lra0.
\]
Then a commutative diagram 
\[
\xymatrix{
H^d_\zar(\cY,\omega_{\cY/W[[t]]}^\bullet)\ar[r]\ar[rd]_D&\frac{dt}{t}\ot
H^d_\zar(\cY,\omega_{\cY/W[[t]]}^\bullet)
\ar[r]&H^d_\zar(\cY,\omega_{\cY}^\bullet)\ar[d]^\Res\\
&H^d_\zar(\cY,\omega_{\cY/W[[t]]}^\bullet)\ar[u]^\cong\ar[r]_-R&\bigoplus_{\ul\nu\in G} W\cdot P_{\ul\nu}
}\]
implies $R\circ D=0$.

We turn to the proof of $C'=1$.
Since $D\wh\eta_{-i_0,\ldots,-i_d}=0$, we have 
$Q(D^i\wh\omega_{i_0\ldots i_d},\wh\eta_{-i_0,\ldots,-i_d})=
D^iQ(\wh\omega_{i_0\ldots i_d},\wh\eta_{-i_0,\ldots,-i_d})=
D^iQ(\omega_{i_0\ldots i_d},\eta_{-i_0,\ldots,-i_d})=0$ by Lemma \ref{unit-lem2} (3).
Therefore the submodule
\[
 VH_{i_0\ldots i_d}(U_K/A_K)_{K((t))^\wedge}\subset \bigoplus_{i=0}^d K((t))^\wedge \cdot D^i\wh\omega_{i_0\ldots i_d}
\] is generated by $\{D^i\wh\omega_{i_0\ldots i_d}\}_{i>0}$, and hence
is annihilated by $R$.
To show $C'=1$, it is enough to show that
\[
R\circ\Phi_{(\cY,D)}(\omega_{p^{-1}i_0\ldots p^{-1}i_d})=p^dR(\omega_{i_0\ldots i_d}).
\]
Let $j_k$ be the unique integer such that $j_k\equiv p^{-1}i_k$ mod $n_k$
and $0<j_k<n_k$. We have
\begin{align*}
R(\omega_{p^{-1}i_0\ldots p^{-1}i_d})
&=\Res\left(n_0^{-1}x_0^{j_0-n_0}x_1^{j_1-1}\cdots x_d^{j_d-1}
\frac{dx_1\cdots dx_d}{(1-x_1^{n_1})\cdots(1-x_d^{n_d})}\frac{dt}{t}\right)\\
&=\frac{(-1)^d}{n_0\cdots n_d}
\left(\sum_{\ul\nu\in G}\nu_0^{j_0}\cdots \nu_d^{j_d}P_{\ul\nu}\right).
\end{align*}
Since $R\circ\Phi_{(\cY,D)}=p^dF\circ R$, we have
\begin{align*}
R\circ\Phi_{(\cY,D)}(\omega_{p^{-1}i_0\ldots p^{-1}i_d})
&=p^dF(R(\omega_{p^{-1}i_0\ldots p^{-1}i_d}))
=p^d\frac{(-1)^d}{n_0\cdots n_d}
\left(\sum_{\ul\nu\in G}\nu_0^{pj_0}\cdots \nu_d^{pj_d}P_{\ul\nu}\right)\\
&
=p^d\frac{(-1)^d}{n_0\cdots n_d}
\left(\sum_{\ul\nu\in G}\nu_0^{i_0}\cdots \nu_d^{i_d}P_{\ul\nu}\right)\\
&
=p^dR(\omega_{i_0\ldots i_d})
\end{align*}
as required.
\end{pf}

\begin{cor}[Unit Root Formula]\label{RC-cor1}
\[
\Phi(\eta_{p^{-1}i_0\ldots p^{-1}i_d})=
\cF_{\check{\ul a}}^{\Dw,\sigma}(t)\ur.
\]
\end{cor}
\begin{pf}
Recall \eqref{RC-lem1-eq2},
\[
\Phi(\eta_{p^{-1}i_0\ldots p^{-1}i_d})
=C
\cF^{\Dw,\sigma}_{\check{\ul a}}(t)
\eta_{i_0\ldots i_d}.
\]
 We want to show $C=1$.
By {\bf(Q5)}, we have
\[Q(\Phi(\omega_{-p^{-1}i_0,\ldots, -p^{-1}i_d}),
\Phi(\eta_{p^{-1}i_0\ldots p^{-1}i_d}))
=p^d\sigma Q(\omega_{-i_0,\ldots,-i_d},\ur)
=p^dQ(\omega_{-i_0,\ldots,-i_d},\ur),
\]
where the second equality follows from the fact $Q(\omega_{-i_0,\ldots,-i_d},\ur)\in\Q^\times$ (Lemma \ref{unit-lem2} (3)).
Applying Theorem \ref{RC-thm1} and \eqref{RC-lem1-eq2} to this, we conclude $C=1$.
\end{pf}
\begin{cor}\label{RC-cor2}
$Q(\eta_{-i_0,\ldots,-i_d},\ur)=0$. In other words, 
$\ur\in VH_{i_0\ldots i_d}(U_K/A_K)_{\wh B_K}$.
\end{cor}
\begin{pf}
Applying Corollary \ref{RC-cor1} to the equality
\[Q(\Phi(\eta_{-p^{-1}i_0,\ldots, -p^{-1}i_d}),
\Phi(\eta_{p^{-1}i_0\ldots p^{-1}i_d}))=p^d\sigma Q(\eta_{-i_0,\ldots,-i_d},\ur)\]
in {\bf(Q5)}, one has
\[
\cF_{{\ul a}}^{\Dw,\sigma}(t)
\cF_{\check{\ul a}}^{\Dw,\sigma}(t)Q(\eta_{-i_0,\ldots,-i_d},\ur)=
p^d\sigma Q(\eta_{-i_0,\ldots,-i_d},\ur).
\]
Comparing the sup norm on the both side, it turns out that $Q(\eta_{-i_0,\ldots,-i_d},\ur)=0$.
\end{pf}
\begin{rem}
There is an alternative proof of Corollary \ref{RC-cor2} with use of the Hodge theory
(monodromy weight filtration).
\end{rem}
\section{Higher Ross symbols and Syntomic regulators}\label{1ext-sect}
In this section we shall discuss the syntomic regulator of the higher Ross symbols. 
The main results are Theorems \ref{main-1} and \ref{main-2}, which are
the $p$-adic counterparts of \cite[Theorem 5.5]{As-Ross1},
and also a generalization of the results in \cite[\S 4.4]{New} in higher
dimension.
\subsection{Higher Ross symbols \cite{As-Ross1}}\label{Rev-Ross-sect}
Let $n_0,\ldots,n_d>1$ be integers.
Let $p$ be a prime number such that
$p\nmid n_0\cdots n_d$.
Let $W=W(\ol\F_p)$ be the Witt ring, and $F$ the $p$-th Frobenius on $W$.
Let $A=W[t,(t-t^2)^{-1}]$, $S=\Spec A$ and 
\begin{equation}\label{HG.5.1}
U=\Spec A[x_0,\ldots,x_d]/((1-x_0^{n_0})\cdots(1-x_d^{n_d})-t)
\end{equation}
the hypergeometric scheme over $A$.
Then,
for $\nu_k\in\mu_{n_k}(A)$ $(k=0,1,\ldots,d$),
the {\it higher Ross symbol} is defined to be
a Milnor symbol
\begin{align}
\Ross:=\left\{\frac{1-x_0}{1-\nu_0 x_0},\ldots,\frac{1-x_d}{1-\nu_d x_d}\right\}
\in K_{d+1}^M(\O(U)),
\end{align} 
in the Milnor $K$-group of $\O(U)$.
We also think it of being an element of Quillen's higher $K$-group
$K_{d+1}(U)$ by the natural map $K^M_i(\O(U))\to K_i(U)$, and
the element is denoted by the same notation.
Let $X\supset U$ be a smooth compactification $X\supset U$ 
such that $Z=X\setminus U$ is a 
relative simple NCD over $A$ (\cite[Proposition 2.1]{As-Ross1}).
Then we expect 
\begin{equation}\label{boundary}
\Ross\in\Image[K_{d+1}(X)^{(d+1)}\to K_{d+1}(U)^{(d+1)}].
\end{equation}
See \cite[4.2]{As-Ross1} for more details.
This is true if $d\leq 2$ (\cite[Corollary 4.4]{As-Ross1}).

\subsection{Category of filtered $F$-isocrystals \cite{AM}}\label{cat-sect}
In \S \ref{regsyn-sect},
we shall employ the category of filtered $F$-isocrystals introduced in \cite[\S 2.1]{AM}
as a fundamental material.
We here recall the notation and some results which we shall need in below.
For a moment, we work over an arbitrary smooth affine variety $S=\Spec(B)$
over $W=W(\ol\F_p)$.
We denote by $B^{\dag}$ the weak completion of $B$.
Namely if $B=W[T_1,\cdots,T_n]/I$, then 
$B^\dag=W[T_1,\cdots,T_n]^\dag/I$ where $W[T_1,\cdots,T_n]^\dag$
is the ring of power series $\sum a_\alpha T^\alpha$ such that
for some $r>1$,
$|a_\alpha|r^{|\alpha|}\to0$ as $|\alpha|\to\infty$. 
Let $K:=\Frac(W)$ be the fractional field, and write $B^\dag_K=K\ot _W B^{\dag}$.

Let $\sigma\colon B^{\dag}\to B^{\dag}$ be a $p$-th Frobenius compatible with
the Frobenius $F$ on $W$.
We define the category $\FilFMIC(S,\sigma)$ (which we call the category of
{\it filtered $F$-isocrystals} on $S$) as follows.
The induced endomorphism $\sigma\otimes_{\bZ}\bQ\colon B_K^{\dag}\to B_K^{\dag}$ 
is also denoted by $\sigma$. 
    An object of $\FilFMIC(S,\sigma)$ 
is a datum $H=(H_{\dR}, H_{\rig}, c, \Phi, \nabla, \Fil^{\bullet})$, where
    \begin{itemize}
            \setlength{\itemsep}{0pt}
        \item $H_{\dR}$ is a coherent $B_K$-module,
        \item $H_{\rig}$ is a coherent $B^{\dag}_K$-module,
        \item $c\colon H_{\dR}\otimes_{B_K}B^{\dag}_K\xrightarrow{\,\,\cong\,\,} H_{\rig}$ is a $B^{\dag}_K$-linear isomorphism,
        \item $\Phi\colon \sigma^{\ast}H_{\rig}\xrightarrow{\,\,\cong\,\,} H_{\rig}$ is an isomorphism of $B^{\dag}_K$-algebra,
        \item $\nabla\colon H_{\dR}\to \Omega_{B_K}^1\otimes H_{\dR}$ is an integrable connection and
        \item $\Fil^{\bullet}$ is a finite descending filtration on $H_{\dR}$ of locally free $B_K$-module (i.e. each graded piece is locally free),
    \end{itemize}
    that satisfies $\nabla(\Fil^i)\subset \Omega^1_{B_K}\ot \Fil^{i-1}$ and
    the compatibility of $\Phi$ and $\nabla$,
namely $\Phi\nabla_\rig=\nabla_\rig\Phi$ where $\nabla_{\rig}\colon H_{\rig}\to\Omega^1_{B_K^{\dag}}\otimes H_{\rig}$ is the connection induced from $\nabla$
under the comparison map $c$.
In what follows we write $\nabla_\rig=\nabla$ for simplicity of notation.

The category $\FilFMIC(S,\sigma)$ is an exact category ({\it not} an abelian category)
in which the tensor product $\ot$ is defined in the customary way.
There are the Tate objects $\O_S(n)=(B_K, B_K^{\dag}, c, p^{-n}\sigma_B, d, \Fil^{\bullet})$, which is the counterpart of the $l$-adic sheaf $\Q_l(n)$. We abbreviate $\O_S=\O_S(0)$.
We write $H(n):=H\ot\O_S(n)$ for an object $H\in \FilFMIC(S,\sigma)$.

Let $u\colon U\to S=\Spec(B)$ be a smooth morphism of smooth $W$-schemes of pure relative dimension.
Assume that there is
a projective smooth morphism $u_X\colon X\to S$ that extends $u$, and
$D:=X\setminus U$ is a normal crossing divisor
with $S$-smooth componentsm (abbreviated relative simple NCD over $S$).
Then the rigid cohomology $H^i_{\rig}(U_{\ol\F_p}/S_{\ol\F_p})$ is defined.
Let
\begin{equation}\label{comparison-eq1}
c:B^\dag_K\ot H^i_\dR(U_K/S_K)\lra H^i_{\rig}(U_{\ol\F_p}/S_{\ol\F_p}).
\end{equation}
be the comparison map (e.g. \cite[2.5]{AM}).
If this is bijective, 
then one can define an object
\[
    H^i(U/S)=(H^i_\dR(U_K/S_K), H^i_{\rig}(U_{\ol\F_p}/S_{\ol\F_p}), c, \nabla, \Phi, \Fil^{\bullet})
\]
of $\Fil$-$F$-$\MIC(S,\sigma)$, where $\nabla$ is the Gauss-Manin
connection, $\Phi$ is the $\sigma$-linear $p$-th Frobenius on the rigid cohomology
and $\Fil^\bullet
H^i_\dR(U_K/S_K)$ is the Hodge filtration.

\medskip

Suppose that $U$ is affine and 
the $i$-th relative rigid cohomology sheaf $R^if_{\rig}j^{\dag}_U\O_{U_K^{\an}}$
is a coherent $j^{\dag}_S\O_{S_K^{\an}}$-module for each $i$.
Then
the comparison map \eqref{comparison-eq1} is bijective for all $i$.
Let $n\geq0$ be an integer.
If $\Fil^{n+1}H^{n+1}_\dR(U_K/S_K)=0$,
then it follows from \cite[Theorem 2.23]{AM} that
we have the symbol map 
\begin{equation}\label{symbol-map}
[-]_{U/S}:K_{n+1}^M(\O(U))\lra \Ext^1_{\FilFMIC(S,\sigma)}(\O_S,H^n(U/S)(n+1))
\end{equation}
to the group of 1-extensions in $\FilFMIC(S,\sigma)$.
\subsection{Syntomic regulators of Higher Ross symbols}\label{regsyn-sect}
We turn to the setting in \S \ref{Rev-Ross-sect}.


\begin{lem}\label{completion-lem}
Suppose $p>\mathrm{lcm}(n_0,\ldots,n_d)$.
Then, for the hypergeometric scheme $U$ in \eqref{HG.5.1},
$R^if_{\rig}j^{\dag}_U\O_{U_K^{\an}}$
is a coherent $j^{\dag}_S\O_{S_K^{\an}}$-module for each $i$, and therefore
the symbol map \eqref{symbol-map} is defined.
\end{lem}
\begin{pf}
Let $n:=\mathrm{lcm}(n_0,\ldots,n_d)$, and let $U'$ be 
the hypergeometric scheme defined by
$(1-x_0^n)\cdots(1-x_d^n)=t$. Then there is a finite etale covering 
$U'\to U$, so that the proof is reduced to the case of $U'$ (e.g. \cite[7.4.1]{CT}).
Therefore, we may assume $n_0=\cdots=n_d$.

We construct a diagram in Theorem \ref{coh-thm} (2) for $U/ A$.
It follows from \cite[Prop. 2.1, Lemmas 2.2, 2.3]{As-Ross1} that
there is a cartesian diagram
\[
\xymatrix{
U\ar[r]^\subset\ar[rd]&X\ar[d]\ar[r]\ar@{}[rd]|{\square}&Y\ar[d]^f\\
&\Spec A\ar[r]^-\subset&\P^1_W
}
\]
such that the following conditions hold.
Let $D_0,D_1,D_\infty$ denote the
scheme-theoretic fiber of $f$ at $t=0$, $t=1$, $t=\infty$ respectively.
\begin{itemize}
\item[\rm(1)]
$Y$ is smooth projective over $W$.
\item[\rm(2)]
$D_0$ is a reduced and relative simple NCD over $W$.
\item[\rm(3)]
$D_\infty$ is a relative simple NCD over $W$ with multiplicity
$\leq n$.
\item[\rm(4)]
$D_1$ is a reduced and irreducible divisor which is 
smooth over $W$ outside the point
$O:=\{x_0=\cdots=x_d=0\}$.
\item[\rm(5)]
Let $\bar Z$ be the closure of $Z:=X\setminus U$ in $Y$.
Then $D_0+(D_1\setminus O)+D_\infty+\bar Z$ is a relative simple NCD over $W$.
\end{itemize}
Let $\bar X\to Y$ be the blowing-up at $O$.
Then 
 the inverse image of $D_1$ has two irreducible components, and
they are a relative simple NCD over $W$ (here we use the assumption $n_0=\cdots=n_d$).
Therefore the diagram
\[
\xymatrix{
U\ar[r]^-\subset\ar[rd]&X\ar[d]\ar[r]&\bar X\ar[d]^{\bar f}\\
&\Spec A\ar[r]&\P^1_W
}
\]
satisfies all conditions in Theorem \ref{coh-thm} (2), and hence the bijectivity
of \eqref{comparison-eq1} follows. This completes the proof.
\end{pf}
By Lemma \ref{completion-lem}, one has the object $H^i(U/S)(r)$ in $\FilFMIC(S,\sigma)$,
and 
the higher Ross symbol $\Ross$ defines
a 1-extension \begin{equation}
0\lra H^d(U/S)(d+1)\lra M_\Ross(U/S)\lra \O_S\lra0
\end{equation}
in $\FilFMIC(S,\sigma)$ by the symbol map \eqref{symbol-map}.
In what follows we take $\sigma$ to be 
the $p$-th Frobenius on $W[[t]]$ given by $\sigma(t)=ct^p$
with some $c\in 1+pW$.
Let $\Phi$ be the Frobenius on $M_\Ross(U_{\ol\F_p}/S_{\ol\F_p})_\rig\cong
M_\Ross(U_K/S_K)_\dR\ot_{A_K}A^\dag_K$ that is $\sigma$-linear.
Let $\Phi_{U/S}$ be the Frobenius on $H^d(U_{\ol\F_p}/S_{\ol\F_p})_\rig$
(without Tate twist).
Notice that
\[
\Phi|_{H^d(U_{\ol\F_p}/S_{\ol\F_p})_\rig}=p^{-d-1}\Phi_{U/S}
\]
by the definition of the Tate twist.
Let
$
e_\Ross\in \Fil^0M_\xi(U/S)_\dR
$
be the unique lifting of $1\in \O(S)$.
Then one has a class $e_\Ross-\Phi(e_\Ross)\in H_\dR^d(U/S)\ot_AA^\dag_K$.
\begin{lem}\label{regsyn-lem1}
$e_\Ross-\Phi(e_\Ross)\in W_dH_\dR^d(U/S)\ot_AA^\dag_K$.
\end{lem}
\begin{pf}
Let $\sigma_i(\ve_i)=(1,\ldots,\ve_i,\ldots,1)\in G=\mu_{n_0}\times\cdots\times\mu_{n_d}$ 
for a fixed primitive $n_i$-th root of
unity $\ve_i\in\mu_{n_i}$.
Put
\[
h:=\prod_{i=0}^d
(1+\sigma_i(\ve_i)+\cdots+\sigma_i(\ve_i)^{n_i-1})\in \Q[G].
\]
Then since
\[
W_dH_\dR^d(U/S)=\ker(h:H_\dR^d(U/S)\to H_\dR^d(U/S)  )\]
by Theorem \ref{dim-cor1} (3), it is enough to show $h(e_\Ross-\Phi(e_\Ross))=0$. 
However since the symbol 
map \eqref{symbol-map} is compatible with respect to the action of
$\Q[G]$ and $\Phi$ commutes with $G$, 
this agrees with $e_{h(\Ross)}-\Phi(e_{h(\Ross)})$.
One can directly show $h(\Ross)=0$ in $K^M_{d+1}(\O(U))\ot\Q$ by defintion,
and hence the vanishing follows.
\end{pf}
\begin{thm}\label{main-1}
Suppose $p>\mathrm{lcm}(n_0,\ldots,n_d)$.
Then  
\[
e_{\Ross}-\Phi(e_\Ross)\equiv 
\sum_{i_0=1}^{n_0-1}\cdots
\sum_{i_d=1}^{n_d-1}(1-\nu_0^{i_0})\cdots(1-\nu_d^{i_d})\cF^{(\sigma)}_{\ul a}(t)\nur
\mod \bigoplus_{i_0,\ldots,i_d}VH_{i_0\ldots i_d}(U/S)_{\wh B_K}
\]
where we put $a_k:=1-i_k/n_k$ and $\ul a:=(a_0,\ldots,a_d)$.
\end{thm}
\begin{pf}
One directly has (cf. \cite[Lemma 4.1]{As-Ross1})
\[
\dlog(\Ross)=(-1)^d\sum_{i_0=1}^{n_0-1}\cdots
\sum_{i_d=1}^{n_d-1}(1-\nu_0^{i_0})\cdots(1-\nu_d^{i_d})
\omega_{i_0\ldots i_d}\frac{dt}{t}.
\]
Let $D:=t\frac{d}{dt}$ be the differential operator acting on 
$M_\xi(U/S)_\dR$.
It follows from \cite[(2.30)]{AM} that we have
\[
D(e_\Ross)=\sum_{i_0=1}^{n_0-1}\cdots
\sum_{i_d=1}^{n_d-1}(1-\nu_0^{i_0})\cdots(1-\nu_d^{i_d})
\omega_{i_0\ldots i_d}.
\]
Let us write
\[
e_\Ross-\Phi(e_\Ross)\equiv \sum_{i_0=1}^{n_0-1}\cdots
\sum_{i_d=1}^{n_d-1}E_{i_0\ldots i_d}(t)\wh\omega_{i_0\ldots i_d}\mod 
VH_{i_0\ldots i_d}(U/S)_{\wh B_K}
\]
where
$\wh\omega_{i_0\ldots i_d}:=F_{\ul a}(t)^{-1}\omega_{i_0\ldots i_d}$.
Then
\begin{align*}
D(e_\Ross-\Phi(e_\Ross))&=D(e_\Ross)-p\Phi D(e_\Ross)\\
&=\sum_{0<i_k<n_k}(1-\nu_0^{i_0})\cdots(1-\nu_d^{i_d})
\bigg(\omega_{i_0\ldots i_d}-p^{-d}\Phi_{U/S}(\omega_{p^{-1}i_0\ldots p^{-1}i_d})\bigg)\\
&\equiv\sum_{0<i_k<n_k}(1-\nu_0^{i_0})\cdots(1-\nu_d^{i_d})
(1-\cF^{\Dw,\sigma}_{\ul a}(t)^{-1})\omega_{i_0\ldots i_d}\\
&\hspace{2cm}\mod \bigoplus VH_{i_0\ldots i_d}(U/S)_{\wh B_K}
\quad \text{(by Theorem \ref{RC-thm1})}\\
&=\sum_{0<i_k<n_k}(1-\nu_0^{i_0})\cdots(1-\nu_d^{i_d})
(F_{\ul a}(t)-F_{\ul a'}(t^{\sigma}))\wh\omega_{i_0\ldots i_d}
\end{align*}
where
$\ul a'=(a'_0,\ldots,a'_d)$ is the Dwork prime.
Hence
we have
\[
\sum_{0<i_k<n_k}\left(t\frac{d}{dt}E_{i_0\ldots i_d}(t)\right)\wh\omega_{i_0\ldots i_d}
=
\sum_{0<i_k<n_k}(1-\nu_0^{i_0})\cdots(1-\nu_d^{i_d})
(F_{\ul a}(t)-F_{\ul a'}(t^{\sigma}))\wh\omega_{i_0\ldots i_d}
\]
as $D\wh\omega_{i_0\ldots i_d}\equiv 0$ mod $VH_{i_0\ldots i_d}(U/S)_{\wh B_K}$.
This implies
\[
E_{i_0\ldots i_d}(t)=C+\int_0^t F_{\ul a}(t)-F_{\ul a'}(t^{\sigma})\frac{dt}{t}
\]
with some constant $C$.
The rest is to show
$C=\psi_p(a_0)+\cdots+\psi_p(a_d)+(d+1)\gamma_p
-p^{-1}\log(c)$.
Notice that $E_{i_0\ldots i_d}(t)$ satisfies that
$E_{i_0\ldots i_d}(t)/F_{\ul a}(t)$ belongs to the ring
$\wh B_K$ (Lemma \ref{unit-lem1}).
If $C=\psi_p(a_0)+\cdots+\psi_p(a_d)+(d+1)\gamma_p
-p^{-1}\log(c)$, then this is satisfied (Theorem \ref{log-thm1}). Suppose that
there is another $C'$ such that 
$E_{i_0\ldots i_d}(t)/F_{\ul a}(t)$ belongs to the ring $\wh B_K$.
Then it follows that so does $(C-C')/F_{\ul a}(t)$ and hence 
$F_{\ul a}(t)\in \wh B_K$.
We show that this is impossible, which completes the proof of Theorem \ref{main-1}.
Let $\sigma_1$ be the Frobenius given by $t^{\sigma_1}=t^p$, and $\Phi_{U/S,\sigma_1}$ the
$\sigma_1$-linear Frobenius on $H^\bullet_\rig(U_\F/S_\F)$.
Let $m>0$ is an integer such that $p^m\equiv 1$ mod $n_k$ for all $k$.
Note $\ul a^{(m)}=\ul a$.
It follows from the unit root formula (Corollary \ref{RC-cor1}) that
we have
\begin{align*}
(\Phi_{U/S,\sigma_1})^m(\eta_{-i_0,\ldots, -i_d})
&=(\Phi_{U/S,\sigma_1})^m(\eta_{-p^{-m}i_0,\ldots,-p^{-m}i_d})\\
&=
\left(\prod_{i=0}^{m-1}
\cF_{\ul a^{(m-i-1)}}^{\Dw,\sigma_1}(t^{p^{m-i-1}})\right)\eta_{-i_0,\ldots, -i_d}\\
&=
\left(\prod_{i=0}^{m-1}
\frac{F_{\ul a^{(m-i-1)}}(t^{p^{m-i-1}})}{F_{\ul a^{(m-i)}}(t^{p^{m-i}})}
\right)\eta_{-i_0,\ldots, -i_d}\\
&=
\left(\frac{F_{\ul a}(t)}{F_{\ul a}(t^{p^m})}
\right)\eta_{-i_0,\ldots, -i_d}.
\end{align*}
Recall $\wh B_K=K\ot_W(W[t,(t-t^2)^{-1},h(t)^{-1}]^\wedge)$ (see \eqref{wh-B} for the
notation).
Choose a $\ell$-th root $\zeta\in W^\times$ of unity with $p\nmid \ell$
such that $(1-\zeta)h(\zeta)\not\equiv 0$ mod $p$.
Let $U_{\F,\zeta}:=U_\F\times_{S_\F} \Spec\F[t]/(t-\zeta)$ be the fiber at $t=\zeta$.
Replacing $m$ with $m\varphi(\ell)$, we may assume $\ell| p^m-1$ and hence $\zeta^{p^m}=\zeta$.
Then the above formula implies that the evaluation
\[
\frac{F_{\ul a}(t)}{F_{\ul a}(t^{p^m})}\bigg|_{t=\zeta}
\]
is the eigenvalue of the $p^m$-th Frobenius
on the rigid cohomology $W_dH^d_\rig(U_{\F,\zeta}/\F)$ with respect to the unit root vector
$\eta_{-i_0,\ldots,-i_d}$. 
Now suppose $F_{\ul a}(t)\in \wh B_K$.
Then
\[
\frac{F_{\ul a}(t)}{F_{\ul a}(t^{p^m})}\bigg|_{t=\zeta}
=\frac
{F_{\ul a}(t)|_{t=\zeta}}
{F_{\ul a}(t^{p^m})|_{t=\zeta}}=1
\]
which contradicts with the Riemann-Weil hypothesis.
This shows $F_{\ul a}(t)\not\in \wh B_K$ as required.
\end{pf}
\begin{rem}\label{remark-1}
The main result of \cite{As-Ross1} is an explicit formula of the pairing
\[
\langle\reg_B(\Ross)\mid\Delta_t\rangle
\]
of the Beilinson regulator $\reg_B(\Ross)$ with a certain homology cycle $\Delta_t$.
See loc. cit. Theorem 5.5 for the details.
One can think
\[
Q(\reg_\syn(\Ross),\whur)
\] 
of being a $p$-adic counterpart of $\langle\reg_B(\Ross)\mid\Delta_t\rangle$
in the following way.
Let $F^\an_{\ul a}(t)$ denote the complex analytic function defined by the 
hypergeometric series $F_{\ul a}(t)$.
Let $y_i^\an$ be the analytic functions defined in the same way as in
Proposition \ref{hat-unit-prop}, and put
\[
\whur^\an:=y^\an_0\nur+y^\an_1D\nur+\cdots+y^\an_dD^d\nur,
\]
and $\ur^\an:=F^\an_{\check{\ul a}}(t)^{-1}\,\whur^\an$.
Proposition \ref{hat-unit-prop} asserts $D\whur^\an=0$.
It follows from Lemma \ref{unit-lem2} (3) that $
Q(\omega_{i_0\ldots i_d},\eta_{-i_0,\ldots,-i_d}^\an)$ is a constant and hence
\[
Q(D^j\omega_{i_0\ldots i_d},\wh\eta_{-i_0,\ldots,-i_d}^\an)=
D^jQ(\omega_{i_0\ldots i_d},\wh\eta_{-i_0,\ldots,-i_d}^\an)=D^jF^\an_{\ul a}(t)
\times\text{(const.)}
\]
for all $j\geq 0$.
On the other hand, the homology cycle $\Delta_\alpha$ satisfies (\cite[Thm.3.1]{As-Ross1})
\[
\langle\omega_{i_0\ldots i_d}\mid\Delta_t\rangle=F^\an_{\ul a}(t)\times
\text{(const.)}
\]
which implies
\[
\langle D^j\omega_{i_0\ldots i_d}\mid\Delta_t\rangle=D^jF^\an_{\ul a}(t)
\times\text{(const.)}
\]
for all $j\geq0$. We thus have
\[
Q(-,\wh\eta_{-i_0,\ldots,-i_d}^\an)=
\langle -\mid\Delta_t\rangle\times\text{(const.)}
\]

\end{rem}
\begin{thm}\label{main-2}
Suppose that $d=1$ and $p\nmid n_0\cdots n_d$ or that $p>\mathrm{lcm}(n_0,\ldots,n_d)$.
Suppose further that $p>d+1$.
Let $\alpha\in W$ satisfy $\alpha\not\equiv 0,1$ and $h(\alpha)\not\equiv 0$ modulo $p$,
so that $(t-\alpha)$ is a maximal ideal of $\wh B_K$.
Let $\sigma$ be the $p$-th Frobenius on $W[[t]]$ given by $\sigma(t)=\alpha^F
\alpha^{-p}t^p$ where $F$ is the $p$-th Frobenius on $W$.
Let 
\[
\reg_\syn:K_{d+1}(U_\alpha)\lra H^{d+1}_\syn(U_\alpha,\Q_p(d+1))\cong 
H^d_\dR(U_{\alpha,K}/K)
\]
be the syntomic regulator map, cf. \eqref{intro-reg}.
Then $\reg_\syn(\Ross|_{U_\alpha})\in W_dH^d_\dR(U_{\alpha,K}/K)$ and 
\[
\frac{Q(\reg_\syn(\Ross|_{U_\alpha}),\eta_{-i_0,\ldots,-i_d})}
{Q(\omega_{i_0\ldots i_d},\eta_{-i_0,\ldots,-i_d})}=
(1-\nu_0^{i_0})\cdots(1-\nu_d^{i_d})\cF^{(\sigma)}_{\ul a}(t)|_{t=\alpha}.
\]
\end{thm}
\begin{pf}
By \cite[Theorem 3.7]{AM},
this is immediate from Theorem \ref{main-1} together with Lemma \ref{regsyn-lem1}. 
\end{pf}

\section{$p$-adic Beilinson conjecture for $K_3$ of K3 surfaces}\label{K3-sect}
We discuss the syntomic regulators for singular K3 surfaces over $\Q$.
Here we mean by {\it singular K3} a (smooth) K3 surface $X$ over a field $k$ 
such that the Picard number of $\ol X$ is maximal.
For a singular K3 surface $X$ over $\Q$, there is a Hecke eigenform $f$ which provides the $L$-function of $X$.
One can formulate the $p$-adic Beilinson conjecture using the $p$-adic $L$-function of $f$.

\medskip

For a smooth projective variety $S$ over a field $k$, we denote by $\NS(S)$ the
Neron-Severi group.
Note $\NS(S)\ot\Q=(\NS(S\times_k \ol k)\ot\Q)^{\mathrm{Gal}(\ol k/k)}$
the fixed part by the Galois group.
Let $h^2_\tr(S,\Q(m))=h^2(S,\Q(m))/\NS(S)\ot\Q(m-1)$ denote the transcendental part of
the motive $h^2(S,\Q(m))$ (cf. \cite[7.2.2]{KMP}).
When $k=\Q$, one can define the $L$-function of $h^2_\tr(S,\Q)$, which we denote by
$L(h^2_\tr(S),s)$.
\subsection{Singular K3 surfaces over $\Q$}
A K3 surface $X$ over a field $k$
is called {\it singular} if the rank of $\NS(X\times_k\ol k)$ is the largest.
We work over the base field $k=\Q$, then $X$ is singular if and only if
the rank is $20$, and hence the transcendental motive $h^2_\tr(X)=h^2_\tr(X,\Q)$ is 2-dimensional.
Every such K3 admits the {\it Shioda-Inose structure} by elliptic curves $E,E'$
with 
complex multiplications (\cite{I-S}, \cite{morrison}),
\begin{equation}\label{Shioda-Inose}
\xymatrix{
X\ar[rd]_{\rho_Z}&&E\times E'\ar[ld]^{\rho_{E\times E'}}\\
&Z
}
\end{equation}
in which the arrows are rational dominant maps of degree 2, and $Z$ is a K3 surface
and the base field is $\ol\Q$.
Moreover $E$ and $E'$ are isogenous.
Thus the $L$-function $L(h^2_\tr(X),s)$ is essentially the $L$-function of
the symmetric square
of the elliptic curves. 
Although the Shioda-Inose structure usually
requires the base field extension, 
the $L$-function has the descent to $\Q$.
\begin{thm}[Livn\'e, \cite{Li}]\label{Li}
Let $X$ be a singular K3 surface over $\Q$.
Then there is a Hecke eigenform $f$ of weight $3$ with complex multiplication
such that $L(h^2_\tr(X),s)=L(f,s)$.
\end{thm}

The $p$-adic Beilinson conjecture asserts that the special values of
$p$-adic $L$-functions are described by syntomic regulators.
See \cite[4.2.2]{Perrin-Riou} and also \cite[Conj.2.7]{Colmez}.
Concerning a singular K3 surface, one can take the $p$-adic $L$-function to be
the $p$-adic $L$-function
\[
L_p(f,\chi,s)
\]
of the associated Hecke eigenform $f$ and some Dirichlet
character $\chi:\Z_p^\times\to\ol\Q^\times_p$ (\cite{MTT}).
See \cite[Conjecture 3.1]{AC} for the general statement of the
$p$-adic Beilinson conjecture for a modular form.
In case of singular K3 surfaces, one can write down
(a part of) the conjecture as follows.

Let $X$ be a K3 surface over $\Q$.
Let $p$ be a prime at which $X$ has a good reduction.
Let $X_{\Z_p}$ be the smooth model over $\Z_p$ and put 
$X_{\Q_p}:=X\times_\Q\Q_p$ and $X_{\F_p}:=X_{\Z_p}\times_{\Z_p}\F_p$.
Let
\[
\reg_\syn:K_3(X_{\Z_p})\lra H^3_\syn(X_{\Z_p},\Q_p(3))\cong H^2_\dR(X_{\Q_p}/\Q_p)
\]
be the syntomic regulator map (\cite{NN}).
Let $X_{\ol\Q_p}:=X_{\Q_p}\times_{\Q_p}\ol\Q_p$ and put
\[
\NS_\dR(X_{\Q_p}):=H^2_\dR(X_{\Q_p}/\Q_p)\cap \NS(X_{\ol\Q_p})\ot\Q_p.
\]
Let $\NS_\dR(X_{\Q_p})^\perp$ denote the orthogonal complement of $\NS_\dR(X_{\Q_p})\subset
H^2_\dR(X_{\Q_p}/\Q_p)$ with respect to the cup-product.
It follows
\[
\NS_\dR(X_{\Q_p})^\perp\op \NS_\dR(X_{\Q_p})\os{\sim}{\lra}H^2_\dR(X_{\Q_p}/\Q_p).
\]
The cup-product induces a non-degenerate pairing on 
\[
\left(H^2_\dR(X_{\Q_p}/\Q_p)/\NS_\dR(X_{\Q_p})\right)\ot\NS_\dR(X_{\Q_p})^\perp
\cong \NS_\dR(X_{\Q_p})^\perp\ot\NS_\dR(X_{\Q_p})^\perp,
\]
which we write by $\langle-,-\rangle$.
\begin{lem}\label{pBconj-lem}
Suppose that $X$ has a good ordinary reduction.
Let
$\alpha_p$ be the eigenvalue of the $p$-th Frobenius on $H^2_\dR(X_{\Q_p}/\Q_p)$
which is a $p$-adic unit (such $\alpha_p$ is unique),
and $\eta\in H^2_\dR(X_{\Q_p}/\Q_p)$ the eigenvector.
Then $\eta\in \NS_\dR(X_{\Q_p})^\perp$.
Let $\omega\in \vg(X_{\Q_p},\Omega^2_{X_{\Q_p}/\Q_p})$ be a non-zero regular $2$-form. 
If $X$ is singular, then $\dim\NS_\dR(X_{\Q_p})^\perp=2$ and hence
\[
\NS_\dR(X_{\Q_p})^\perp=\Q_p\omega+\Q_p\eta.
\]
\end{lem}
\begin{pf}
We want to show that the cup-product $\eta\cup z$ vanishes for any 
$z\in \NS_\dR(X_{\Q_p})$.
The eigenvalues of $\Phi$
on $\NS_\dR(X_{\Q_p})$ are $p\times$(roots of unity).
Take $m>0$ such that $\Phi^m=p^m$ on $\NS(X_{\Q_p})$.
Since $\Phi(x)\cup \Phi(y)=\Phi(x\cup y)=p^2 x\cup y$, we have
\[
\Phi^m(\eta)\cup\Phi^m(z)=p^{2m}(\eta\cup z)\quad\Longleftrightarrow\quad
\alpha_p^mp^m(\eta\cup z)=p^{2m}(\eta\cup z).
\]
Since $\alpha_p\in \Z_p^\times$, it follows that $\eta\cup z=0$.
\end{pf}

\begin{conj}[weak $p$-adic Beilinson conjecture for ordinary singular K3 surfaces]\label{pBconj}
Let $X$ be a singular K3 surface over $\Q$
and $f$ the corresponding Hecke eigenform of weight $3$.
Fix $\omega\in \vg(X,\Omega^2_{X/\Q})$ a regular 2-form.
Suppose that $X$ has a good ordinary reduction at $p$, which is equivalent to
that $p\nmid a_p$ where $f=\sum_{n=1}^\infty a_n q^n$.
Let $\alpha_p$ be the unit root of $T^2-a_pT+p^2$, and 
$\eta\in \NS(X_{\Q_p})$ its eigenvector
(note that $\omega$ is unique up to $\Q^\times$ and $\eta$ is unique
up to $\Q^\times_p$).
Then there is an integral element $\xi\in K_3(X)^{(3)}_\Z$ (cf. \cite{Scholl}) and
a constant $C\in \Q^\times$ not depending on $p$ such that
\begin{equation}\label{pBconj-eq}
L_p(f,\omega_{\mathrm{Tei}}^{-1},0)=C(1-p^2\alpha^{-1}_p)
\frac{\langle\reg_\syn(\xi),\eta\rangle}{\langle\omega,\eta\rangle}
\end{equation}
where $\omega_{\mathrm{Tei}}$ is the Teichm\"uller character.
\end{conj}
\begin{rem}
For a singular K3 surface $X$,
one can expect $K_i(X)^{(j)}=K_i(X)^{(j)}_\Z$ with $i\ne1$ thanks to the
Shioda-Inose structure \eqref{Shioda-Inose}.
See \cite[Remark 6.8]{As-Ross1} for details.
\end{rem}

\subsection{Review of \cite{ono}}\label{ono-sect}
In \cite{ono}, Ahlgren, Ono and Penniston study the $L$-function of 
a K3 surface $Y_a$ over $\Q$ defined by an affine equation
\[
w^2=u_1u_2(1+u_1)(1+u_2)(u_1-a u_2),\quad a\in \Q\setminus\{0,1\}.
\] 
Their first main result is the following.
\begin{thm}[{\cite[Theorem 1.1]{ono}}]\label{d=2-thm1}
Let
\[
E_a:y^2=x\left(x^2+2x-\frac{a}{1-a}\right)
\]
\[
E'_a:(1-a)y^2
=x\left(x^2+2x-\frac{a}{1-a}\right)
\]
be ellptic curves over $\Q$.
Then $L(h^2_\tr(Y_a),s)=L(h^2_\tr(E_a\times E'_a),s)$.
\end{thm}
The explicit correspondence between $Y_a$ and $E_a\times E'_a$
which gives an isomorphism
$h^2_\tr(Y_a)\cong h^2_\tr(E_a\times E'_a)$ of motives over $\Q$
is constructed by van Geemen and Top \cite[Theorem 1.2]{GT}.

\medskip

Following \cite{ono}, we call $Y_a$ {\it modular} if $L(h^2_\tr(Y_a),s)$ is
the $L$-function of a Hecke eigenform of weight $3$ with complex multiplication,
or equivalently $E_a$ has a complex multiplication (\cite[Theorem 1.2]{GT}).
The second main result of \cite{ono} 
gives the complete list of $a$'s for $Y_a$ to be modular.
\begin{thm}[{\cite[Theorem 1.2]{ono}}]\label{d=2-thm2}
The K3 surface $Y_a$ is modular if and only if
$a=-1,4^{\pm1}$,
$-8^{\pm1},64^{\pm1}$.
In each case, $E_a$ has complex multiplication.
\end{thm}
The corresponding Hecke eigenforms are as follows (\cite[p.366--367]{ono}).
Let $\eta(z)$ be the Dedekind eta function. Let
\[
A=\eta^6(4z),\quad B=\eta^2(z)\eta(2z)\eta(4z)\eta^2(8z),\quad
C=\eta^3(2z)\eta^3(6z),\quad D=\eta^3(z)\eta^3(7z)
\]
be weight $3$ newforms of level $16$, $8$, $12$, $7$ respectively.
Let $\chi_D$ denote the quadratic character associated to the quadratic field $\Q(\sqrt{D})$.
Then the corresponding Hecke eigenforms are given as follows.
\begin{equation}\label{d=2-table2}
\text{
 \begin{tabular}{c|ccccccc}
$a$&$-1$&$4$&$1/4$&$-8$&$-1/8$&$64$&$1/64$\\
\hline
Hecke eigenform&$B\ot\chi_{-4}$&$C$&$C\ot\chi_{-4}$&$A$&$A\ot\chi_8$&$D$&$D\ot\chi_{-4}$
\end{tabular}
}
\end{equation}
where $f\ot\chi$ denotes the $\chi$-twist of the modular form.

\subsection{$p$-adic Beilinson conjecture for K3 surfaces in \cite{ono}}\label{pB-sect}
Let
\[
U_a:(1-x_0^2)(1-x_1^2)(1-x_2^2)=a
\]
the hypergeometric scheme over $\Q$, 
and let $X_a$ be a smooth compactification of $U_a$ which
is a K3 surface. Put $Z_a:=X_a\setminus U_a$.
We discuss the $p$-adic Beilinson conjecture for $X_a$.
Let $p>3$ be a prime at which $X_a$ has a good ordinary reduction.
We take integral models $X_{a,\Z_{(p)}}\supset U_{a,\Z_{(p)}}$ which are smooth
over $\Z_{(p)}$ (Proposition \ref{compact-prop}). 
For a $\Z_{(p)}$-ring $R$, we write $X_{a,R}:=X_{a,\Z_{(p)}}\times_{\Z_{(p)}}R$ etc.
Recall the 2-forms
\[
\omega_{1,1,1},\, \eta_{1,1,1}\in W_2H^2_\dR(U_a/\Q_p)\cong
H^2_\dR(X_a/\Q_p)/H^2_{\dR,Z_a}(X_a/\Q_p)
\]
from \eqref{form-eq1} and \eqref{unit-eq3}.
Let
\[
\Ross=\left\{\frac{1-x_0}{1+x_0},\frac{1-x_1}{1+x_1},\frac{1-x_2}{1+x_2}
\right\}\in K_3^M(\O(U_{a,\Z_{(p)}}))
\]
be the higher Ross symbol.
We think $\Ross$ to be an element of $K_3(U_{a,\Z_{(p)}})^{(3)}$ under the natural map
$K_3^M(\O(U_{a,\Z_{(p)}}))\to K_3(U_{a,\Z_{(p)}})^{(3)}$.
There is the exact sequence
\[
K_3^{Z_{a,\Z_{(p)}}}(X_{a,\Z_{(p)}})^{(3)}\lra K_3(X_{a,\Z_{(p)}})^{(3)}
\lra K_3(U_{a,\Z_{(p)}})^{(3)}.
\]
Then $\Ross$ lies in the image of $K_3(X_{a,\Z_{(p)}})^{(3)}$ (\cite[Corollary 4.4]{As-Ross1}),
so that there is a lifting
 $\wt\Ross\in K_3(X_{a,\Z_{(p)}})^{(3)}$.
It is a standard argument on the syntomic regulator maps
to see that there is a commutative diagram
\begin{equation}\label{CD-Z}
\xymatrix{
K_3(X_{a,\Z_{(p)}})^{(3)}\ar[r]^-{\reg_\syn}
\ar[d]& H^3_\syn(X_{a,\Z_p},\Q_p(3))\cong H^2_\dR(X_{a,\Q_p}/\Q_p)\ar[d]\\
K_3(X_{a,\Z_{(p)}})^{(3)}/K_3^{Z_{a,\Z_{(p)}}}(X_{a,\Z_{(p)}})^{(3)}\ar[r] &
H^2_\dR(X_{a,\Q_p}/\Q_p)/\NS_\dR(X_{a,\Q_p})\ar[r]^-{\langle-,\eta_{1,1,1}\rangle}&\Q_p
}\end{equation}
where $\langle-,\eta_{1,1,1}\rangle$ is the cup-product pairing 
which is well-defined by Lemma \ref{pBconj-lem}.
The diagram implies that
$\langle\reg_\syn(\wt\Ross),\eta_{1,1,1}\rangle$ does not depend
on the choice of the lifting.
Applying Theorem \ref{main-2}, we have the description of
the right hand side of \eqref{pBconj-eq} in Conjecture \ref{pBconj} in terms of
our $p$-adic function $\cF^{(\sigma)}_{\ul a}(t)$.
\begin{thm}\label{pB-thm}
Let $\sigma(t)=a^{1-p}t^p$. Then
\[
\frac{\langle\reg_\syn(\wt\Ross),\eta_{1,1,1}\rangle}{\langle\omega_{1,1,1},\eta_{1,1,1}\rangle}
=8\cF^{(\sigma)}_{\frac12,\frac12,\frac12}(t)|_{t=a}.
\]
\end{thm}
Next, we see the left hand side of \eqref{pBconj-eq} in Conjecture \ref{pBconj},
namely the $p$-adic $L$-function.
Recall the K3 surface $Y_a$ from \S \ref{ono-sect}.
There is a dominant rational map 
\begin{equation}\label{XY}
\rho:X_a\lra Y_a
\end{equation}
over $\Q$ (\cite[(6.7)]{As-Ross1}), which induces an isomorphism
\begin{equation}\label{XY-L}
h^2_\tr(Y_a,\Q)\cong
h^2_\tr(X_a,\Q)
\end{equation}
of motives over $\Q$.
Therefore the $p$-adic Beilinson conjecture for $Y_a$ is equivalent to that for $X_a$.
The $L$-functions of $X_a$ and $Y_a$ agree, and if  
$a=-1,4^{\pm1},-8^{\pm1},64^{\pm1}$, then they are the $L$-functions
of the Hecke eigenforms as in the table \eqref{d=2-table2}.
Together with Theorem \ref{pB-thm},
Conjecture \ref{pBconj} for $K_3(X_a)$ or $K_3(Y_a)$ can be formulated as follows.
\begin{conj}
\label{pRZconj}
Let $a=-1,4^{\pm1},-8^{\pm1},64^{\pm1}$.
Let $f_a$ be the Hecke eigenform corresponding to $Y_a$, 
cf. table \eqref{d=2-table2}.
Let $p>3$ be a prime such that $p\nmid a_p$ 
and let $\alpha_p$ be the unit root of $T^2-a_pT+p^2$.
Let $\sigma$ be the $p$-th Frobenius given by $\sigma(t)=a^{1-p}t^p$.
Then there is
a constant $C_a\in \Q^\times$ not depending on $p$ such that
\[
L_p(f_a,\omega_{\mathrm{Tei}}^{-1},0)=C_a(1-p^2\alpha^{-1}_p)
\cF^{(\sigma)}_{\frac12,\frac12,\frac12}(t)|_{t=a}.
\]
\end{conj}
%
%

In view of \cite[Theorem 6.9]{As-Ross1}, it is also plausible to expect
that Conjecture \ref{pRZconj} remains true for $a=1$.
\begin{conj}\label{pRZconj-1}
Let $p\equiv1 $ mod $4$ be a prime.
Let $\alpha_p\in \Z_p$ be the root of $T^2-a_pT+p^2$ such that
$\alpha_p\equiv a_p$ mod $p$ where $A=\eta^6(4z)=\sum a_nq^n$.
Then there is
a constant $C_1\in \Q^\times$ not depending on $p$ such that
\[
L_p(A,\omega_{\mathrm{Tei}}^{-1},0)=C_1(1-p^2\alpha^{-1}_p)
\cF^{(\sigma)}_{\frac12,\frac12,\frac12}(t)|_{t=1}
\]
where $\sigma(t)=t^p$.
\end{conj}

\subsection{Some other elliptic K3 surfaces and $p$-adic regulators}\label{K3-1-sect}
Let $A=\Q[s,(s-s^2)^{-1}]$ and
\[
V_n=\Spec A[x_0,x_1]/((1-x_0^n)(1-x_1^n)-s)
\]
a hypergeometric scheme for $n\geq 2$ an integer.
Let $C_n\supset V_n$ be a smooth compactification over $A$
(Proposition \ref{compact-prop}).
The relative dimension of $C_n/A$ is $1$.
Let $J(C_n)\to\Spec A$ be the jacobian scheme.
Put $H:=H^1_\dR(C_n/A)$.
For $0<i_0,i_1<n$, we denote by $H_{\ol\Q}(i_0,i_1)\subset H_{\ol\Q}:=H\ot_\Q\ol\Q$ the eigenspace 
defined in \eqref{rev-eigen-eq}.
Then the subspace
\[
\sum_{\gcd(r,n)=1}H_{\ol\Q}(ri_0,ri_1)\subset H_{\ol\Q}
\]
is endowed with $A$-module structure, which we denote by $H_A(i_0,i_1)$,
\[
H_A(i_0,i_1)\ot_\Q\ol\Q=\sum_{\gcd(r,n)=1}H_{\ol\Q}(ri_0,ri_1).
\]
Note $H_A(i_0,i_1)=H_A(i'_0,i'_1)$ if and only if $i_0\equiv ri'_0$ and $i_1\equiv ri'_1$ mod $n$
for some $r$ prime to $n$.
Let
\[
H=\bigoplus_{i_0,i_1} H_A(i_0,i_1)
\]
be the decomposition
where $(i_0,i_1)$ runs over representatives of the set
$\{(\bar i_0,\bar i_1)\in (\Z/n\Z)^2\mid\bar i_0\ne0,\,\bar i_1\ne0\}/\sim$
with $(\bar i_0,\bar i_1)\sim(\bar i'_0,\bar i'_1)$
$\Leftrightarrow$ $(\bar i_0,\bar i_1)=(r\bar i'_0,r\bar i'_1)$ for some $r\in(\Z/n\Z)^\times$.
Thanks to the Poincare reducibility theorem (\cite[\S 19 Thm. 1]{Mum-abel}), 
the above decomposition induces an isogeny
\[
J(C_n)\lra\prod_{i_0,i_1}J_{i_0,i_1}
\]
of abelian $A$-schemes.
We note that
\[
\sum_{0<r<n}H_{\ol\Q}(ri_0,ri_1)
=\sum_{0<r<n}H_A(ri_0,ri_1)\ot\ol\Q
\] is associated to 
the hypergeometric curve $y^n=x^{n-i_0}(1-x)^{n-i_1}(1-(1-t)x)^{i_1}$ of Gauss type
(cf. \cite[4.6]{New}).
Hereafter we consider the cases
\[
n=3,4,6.
\]
Let $J_n$ be the associated abelian scheme to $H_A(1,n-1)$.
Since
\begin{equation}\label{K3-1-eq1}
H_A(1,n-1)\ot_\Q\ol\Q=H_{\ol\Q}(1,n-1)\op H_{\ol\Q}(n-1,1)
\end{equation}
is of rank 4 (as $n=3,4,6$), 
the relative dimension of $J_n/A$ is $2$.
Let $\tau$ be the involution on $V_n$ given by $(x_0,x_1)\mapsto (x_1,x_0)$.
This acts on \eqref{K3-1-eq1} such that $\tau(H_{\ol\Q}(1,n-1))=H_{\ol\Q}(n-1,1)$.
Let $H_A(1,n-1)^\pm$ denote the eigenspace of $\tau$ with eigen value $\pm1$,
and $E^\pm_n$ the corresponding decomposition of $E_n$.
Both of $E_n^\pm$ are elliptic curves over $A$.
We thus have surjective morphisms
\begin{equation}\label{K3-1-eq2}
\xymatrix{
F_n^\pm:C_n\ar[r]&E_n^\pm
}
\end{equation}
over $A$, which satisifes
\begin{equation}\label{K3-1-eq3}
\xymatrix{
H^1_\dR(E_n^\pm/A)
\ar[r]^-{(F^\pm_n)^*}_-\cong&
 H_A(1,n-1)^\pm\subset H^1_\dR(C_n/A),
}\end{equation}
\begin{equation}\label{K3-1-eq4}
(i_0,i_1)\ne(1,n-1),(n-1,1)\quad\Longrightarrow\quad
(F_n^\pm)_*(H_{\ol\Q}(i_0,i_1))=0.
\end{equation}
%
Let $B=\Q[t,(t-t^2)^{-1}]$, and  
\[
U_n=\Spec B[x_0,x_1,x_2]/((1-x_0^n)(1-x_1^n)(1-x_2^2)-t).
\]
Let $X_n\supset U_n$ be a smooth compactification over $B$.
Let $S:=\Spec B[x_2,(1-x_2^2)^{-1},(1-t-x_2^2)^{-1}]$.
Let $\phi:S\to\Spec A$ be the morphism given by 
$\phi^*(s)=t(1-x_2^2)^{-1}$, and
\[
\xymatrix{
V_{n,S}\subset C_{n,S}\ar[rd]\ar[rr]^-{F^\pm_{n,S}}&&E_{n,S}^\pm
\ar[ld]^-{h^\pm_{n,S}}\\
&S
}
\]
the base change
of \eqref{K3-1-eq2} by $\phi$,
where 
\begin{align}
V_{n,S}&:=V_n\times_{A}S\label{K3-1-VnS}\\
&=\Spec B[x_0,x_1,x_2,(1-x_2^2)^{-1},(1-t-x_2^2)^{-1}]/(
(1-x_0^n)(1-x_1^n)-t(1-x^2_2)^{-1})\notag\\
&=\Spec B[x_0,x_1,x_2,(1-t-x_2^2)^{-1}]/(
(1-x_0^n)(1-x_1^n)(1-x^2_2)-t)\notag\\
&\hra U_n.\notag
\end{align}
The morphisms $h^\pm_{n,S}$ give rise to
relatively minimal elliptic fibrations
$\cE^\pm_n\to \P^1_B(x_2)$ over $B$,
so that we have a commutative diagram
\begin{equation}\label{K3-1-eq5}
\xymatrix{
V_{n,S} \ar[rd]_{g_n}\ar[rr]^-{f_n^\pm}&&\cE_n^\pm\ar[ld]^{h^\pm_n}\\
&\P^1_B(x_2)
}
\end{equation}
of smooth $B$-schemes.
\begin{prop}\label{ellK3-prop1}
$\cE_n^\pm$ $(n=3,4,6)$ are elliptic K3 surfaces over $B$.
\end{prop}
\begin{pf}
Let $x\in \Spec B$ be a closed point, and $k(x)$ the residue field.
We denote by $\cE_{n,x}^\pm$ the fiber at $x$.
We first show that the proposition is reduced to show that
$\cE_{n,x}^\pm$ are K3 surfaces over $k(x)$ for any $x$.
Indeed, suppose that they are true, namely
$H^1(\cE_{n,x}^\pm,\O)=0$ and $\Omega^2_{\cE_{n,x}^\pm/k(x)}\cong\O$.
Then it follows from \cite[III, 12.9]{Ha} that one has the vanishing
$H^1(\cE_{n}^\pm,\O)=0$.
Moreover $H^0(\Omega^2_{\cE_{n}^\pm/B})$ is locally free 
$B$-module of rank one
(and hence $H^0(\Omega^2_{\cE_{n}^\pm/B})\cong B$), and
the natural map
\[
H^0(\Omega^2_{\cE_{n}^\pm/B})\ot_Bk(x)\lra H^0(\Omega^2_{\cE_{n,x}^\pm/k(x)})
\]
is bijective. 
We have a commutative diagram
\[
\xymatrix{
(H^0(\Omega^2_{\cE_{n}^\pm/B})\ot_B\O_{\cE_n^\pm})\ot_Bk(x)
\ar[d]_\cong\ar[r]&\Omega^2_{\cE_{n}^\pm/B}\ot_Bk(x)\ar[d]^\cong\\
H^0(\Omega^2_{\cE_{n,x}^\pm/k(x)})\ot_{k(x)}\O_{\cE_{n,x}^\pm}\ar[r]&\Omega^2_{\cE_{n,x}^\pm/k(x)}
}
\]
and the bottom arrow is bijective as $\Omega^2_{\cE_{n,x}^\pm/k(x)}\cong\O$.
Hence this implies an isomorphism
\[
\xymatrix{
\O_{\cE_n^\pm}\cong
H^0(\Omega^2_{\cE_{n}^\pm/B})\ot_B\O_{\cE_n^\pm}
\ar[r]^-\cong&\Omega^2_{\cE_{n}^\pm/B}
}
\]
as required.

We may now replace the base ring $B$ with a field $k$ of characteristic zero.
Then we want to show the vanishing $H^1(\cE_n^\pm,\O)=0$
and an isomorphism
$\Omega^2_{\cE_n^\pm/k}\cong \O$.
To do this, we may further replace $k$ with $\C$.
There is the injective map
$(f^-_n)^*:H^1_\dR(\cE_n^\pm)\hra H^1_\dR(X_n)=W_1H^1_\dR(U_n)$.
However since $W_1H^1_\dR(U_n)=0$ (Theorem \ref{dim-cor1} (2)), it turns out
that $H^1_\dR(\cE_n^\pm)=0$ and hence the vanishing 
$H^1(\cE_n^\pm,\O)=0$ follows.
The rest is to show an isomorphism
$\Omega^2_{\cE_n^\pm/\C}\cong \O$.
To do this, we employ the canonical bundle formula.
There is a $A$-rational point of $C_n/A$, and hence 
there is a section of the elliptic fibration $h^\pm_n:\cE^\pm_n\to \P^1_k$.
Hence it follows from the canonical bundle formula \cite[V, (12.3)]{barth} that
one has
\begin{equation}\label{ellK3-prop1-eq1}
K_{\cE^\pm_n}:=\Omega^2_{\cE_n^\pm/k}\cong (h^\pm_n)^*\O_{\P^1}(e)
\end{equation}
with $e=\chi(\O_{\cE^\pm_n})-2\chi(\O_{\P^1})=\chi(\O_{\cE^\pm_n})-2$.
We show $e=0$. 
We denote the topological Euler number of $X$ by $e(X)$.
\begin{align*}
\chi(\O_{\cE^\pm_n})
&=\frac{1}{12}(K_{\cE^\pm_n}^2+e(\cE^\pm_n))&\text{(Noether's
formula, \cite[I, (5.5)]{barth})} \\
&=\frac{1}{12}e(\cE^\pm_n) &\text{(by \eqref{ellK3-prop1-eq1})}\\
&=\frac{1}{12}\sum_{s}e(Z_{n,s}^\pm) &\text{(\cite[III, (11.4)]{barth})}\\
\end{align*}
where $Z_{n,s}^\pm$ runs over all singular fibers of $h^\pm_n:\cE_n^\pm\to\P^1(x_2)$.
To compute the last term, we write up the singular fibers for each $n=3,4,6$.
The singlar fibers of $h_n^\pm$ appear at $x_2=\pm 1,\pm\sqrt{1-t}$ and $x_2=\infty$
(we think $t\in\C\setminus\{0,1\}$ of being a constant).
The Kodaira type of each singular fiber is determined by the local monodromy on
$R^1(h^\pm_n)_*\Q$, and it can be computed from the local monodromy
for the fibration $E_n^\pm\to \Spec A=\bA^1(s)\setminus\{0,1\}$ (here 
we think $s$ of being a parameter).
By virtue of \eqref{K3-1-eq3}, this is isomorphic to the monodromy of
the Gaussian hypergeometric function ${}_2F_1\left({\frac{1}{n},1-\frac{1}{n}\atop 1};s\right)$,
which is well-understood. 
In this way, we obtain the complete list of singular fibers of $h_n^\pm$,
\begin{center}
\begin{tabular}{c||c|c|c}
&$x_2=\pm1$&$x_2=\pm\sqrt{1-t}$&$x_2=\infty$\\
\hline
\hline
$h_3^\pm$&IV$^*$&I$_1$&I$_6$\\
\hline
$h_4^\pm$&III$^*$&I$_2$&I$_2$\\
\hline
$h_6^\pm$&II$^*$&I$_1$&I$_2$\\
\end{tabular}
\end{center}
It is now immediate to have $\sum_{s}e(Z_{n,s}^\pm)=24$ in all cases. Hence we have 
$e=0$ as required.
\end{pf}
\begin{lem}\label{UV-lem}
Let $V_{n,S}=V_n\times_{A}S\hra U_n$ be as in \eqref{K3-1-VnS}. One has
\[
W_2H^2_\dR(U_n/B)=W_2H^2_\dR(V_{n,S}/B).
\]
Hence the pull-back
$(f^\pm_n)^*:H^2_\dR(\cE_n^\pm/B)\lra W_2H^2_\dR(U_n/B)$
is defined.
\end{lem}
\begin{pf}
(This is implicitly proven in the proof of \cite[Theorem 3.4]{As-Ross1}).
It is enough to show the lemma at each fiber, so that we may replace $B$ with
$\C$ and may assume that $U_n$ and $V_{n,S}$ are smooth complex affine varities.
Let $F_1,F_2\subset U_n$ be the complement of $V_{n,S}$ which are irreducible
affine curves
with unique singular points $P_1\in F_1$ and $P_2\in F_2$. 
Put $U_n^\circ:=U_n\setminus\{P_1,P_2\}$ and $F^\circ_i:=F_i\setminus\{P_i\}$. Then
there is an exact sequence
\[
\xymatrix{
&&H^2(U_n)\ar[d]^\cong\\
H^1(V_{n,S})\ar[r]^-{\Res_1}&H^0(F_i^\circ)^\op\ot\Q(-1)\ar[r]&H^2(U_n^\circ)\ar[r]&H^2(V_{n,S})
\ar[r]^-{\Res_2}&H^1(F_i^\circ)^\op\ot\Q(-1)
}
\]
of mixed Hodge structures, where $H^i(-)=H^i(-,\Q)$ denotes the Betti cohomology.
Since $\Res_1$ is surjective and the weight of $H^1(F^\circ_i)\ot\Q(-1)$ is $\geq3$,
we have $W_2H^2(U_n)=W_2H^2(V_{n,S})$ as required.
\end{pf}
Put $H(U_n):=W_2H^2_\dR(U_n/B)$.
The sum 
\[
\sum_rH(U_n)_{\ol\Q}(ri_0,ri_1,1)\subset
H(U_n)_{\ol\Q}:=H(U_n)\ot_\Q\ol\Q
\]
of the subspaces
is endowed with $B$-module structure, which we denote by $H(U_n)_B(i_0,i_1,1)$,
\[
H(U_n)_B(i_0,i_1,1)\ot\ol\Q=
\sum_rH(U_n)_{\ol\Q}(ri_0,ri_1,1),
\]
\[
H(U_n)=\bigoplus_{i_0,i_1}H(U_n)_B(i_0,i_1,1).
\]
Moreover we denote by $H(U_n)_B(i_0,i_1,1)^+$ (resp. $H(U_n)_B(i_0,i_1,1)^-$)
the fixed part (resp. anti-fixed part) by the involution $\tau(x_0,x_1,x_2)=(x_1,x_0,x_2)$.
\begin{prop}\label{ellK3-prop2}
The image of 
\[
(f^\pm_n)^*:H^2_\dR(\cE_n^\pm/B)\lra H(U_n)=W_2H^2_\dR(U_n/B)
\]
agrees with the component $H(U_n)_B(1,n-1,1)^\pm$.
\end{prop}
\begin{pf}
There is a commutative diagram
\begin{equation}\label{K3-1-eq6.0}
\xymatrix{
H^2_\dR(\cE_{n}^\pm/B)\ar[r]^-{(f_n^\pm)^*}\ar[d]&W_2H^2_\dR(V_{n,S}/B)
\os{(\star)}{=}H(U_n)\ar[d]^\cap\\
H^2_\dR(E_{n,S}^\pm/B)\ar[r]^-{(F_{n,S}^\pm)^*} &
H^2_\dR(V_{n,S}/B)\\
H^1_\dR(S,R^1(h_{n}^\pm)_*\Omega^\bullet_{E_{n,S}^\pm/S})\ar[u]^i
\ar[r]& H^1_\dR(S,R^1g_{n,*}\Omega^\bullet_{V_{n,S}/S})\ar[u]_\cong
}
\end{equation}
where the bottom arrow is 
induced from $F_n^\pm$ in \eqref{K3-1-eq2} and $(\star)$ 
follows from Lemma \ref{UV-lem}.
The map $i$ is injective and the cokernel of $i$ is generated by
the cycle class $e$ of a section of $E^\pm_{n,S}\to S$.
One easily sees $(F^\pm_{n,S})^*(e)=0$.
Therefore the image of $F_{n,S}^*$ agrees with
the image of $H^1_\dR(S,R^1(h_{n}^\pm)_*\Omega^\bullet_{E_{n,S}^\pm/S})$, and hence
we have
\[
(f^\pm_n)^*H^2_\dR(\cE_n^\pm/B)\subset H(U_n)_B(1,n-1,1)^\pm.
\]
Since $H(U_n)_B(1,n-1,1)^\pm$ is an irreducible connection (Theorem \ref{dim-cor2}),
the equality holds in the above, as required.
\end{pf}
\begin{cor}
For a geometric point $\bar x\to\Spec B$, let $\cE^\pm_{n,\bar x}$ be the fiber at $\bar x$.
Then the rank of the Neron-Severi group $\NS(\cE^\pm_{n,\bar x})$ is $\geq 19$.
\end{cor}
\begin{pf}
This follows from Proposition \ref{ellK3-prop2} and the fact
$\mathrm{rank}H(U_n)_B(1,n-1,1)^\pm=3$.
\end{pf}
Taking the Hodge filtration $F^2$ of the diagram
\eqref{K3-1-eq6.0}, we have a commutative diagram
\begin{equation}\label{K3-1-eq6}
\xymatrix{
H^0(\cE_n^\pm,\Omega^2_{\cE_n^\pm/B})
\ar@<0.5ex>[r]^-{(f_n^\pm)^*} &
H^0(X_n,\Omega^2_{X_n/B})\ar@<0.5ex>[l]^-{(f_n^\pm)_*}\\
&F^2H^1_\dR(S,R^1g_{n,*}\Omega^\bullet_{V_{n,S}/S})\ar[u]_\cong\\
F^2H^1_\dR(S,R^1(h_{n}^\pm)_*\Omega^\bullet_{E_{n,S}^\pm/S})\ar[uu]^\cong
\ar@<0.5ex>[r]& F^2H^1_\dR(S,R^1g_{n,*}\Omega^\bullet_{C_{n,S}/S})\ar[u]_\cong
\ar@<0.5ex>[l]
}
\end{equation}
together with the push-forward maps.

\begin{lem}\label{ellK3-lem}
Let $\omega_{i_0,i_1,1}\in H^0(X_n,\Omega^2_{X_n/B})$ br the regular 2-forms
in \eqref{form-eq1}. Then
\begin{equation}\label{ellK3-prop-1}
(f_n^\pm)^*H^0(\cE_n^\pm,\Omega^2_{\cE_n^\pm/B})=
B(\omega_{1,n-1,1}\mp\omega_{n-1,1,1})
\end{equation}
\begin{equation}\label{ellK3-prop-2}
(i_0,i_1)\ne(1,n-1),(n-1,1)\quad\Longrightarrow\quad
(f_n^\pm)_*(\omega_{i_0,i_1,1})=0,
\end{equation}
\begin{equation}\label{ellK3-prop-3}
(f_n^\pm)_*(\omega_{1,n-1,1})=\mp
(f_n^\pm)_*(\omega_{n-1,1,1})\ne0.
\end{equation}
\end{lem}
\begin{pf}
Since the involution $\tau(x_0,x_1,x_2)=(x_1,x_0,x_2)$ on $U_n$ satisfies
$\tau\omega_{1,n-1,1}=-\omega_{n-1,1,1}$, one has
\[
F^2H(U_n)_B(1,n-1,1)^\pm=B(\omega_{1,n-1,1}\mp\omega_{n-1,1,1}).
\]
Now \eqref{ellK3-prop-1} is immediate from Proposition \ref{ellK3-prop2}.
Moreover \eqref{ellK3-prop-2} follows from \eqref{K3-1-eq4} by virtue of the diagram 
\eqref{K3-1-eq6}.
By the construction of $F_n^\pm$, they satisfy $(F^+_n)_*(F^-_n)^*=(F^-_n)_*(F^+_n)^*=0$,
which implies $(f^+_n)_*(f^-_n)^*=(f^-_n)_*(f^+_n)^*=0$.
Therefore,
\[
(f_n^+)_*(\omega_{1,n-1,1}+\omega_{n-1,1,1})=0,\quad
(f_n^-)_*(\omega_{1,n-1,1}-\omega_{n-1,1,1})=0.\]
There remains to show the non-vanishing $(f_n^\pm)_*(\omega_{1,n-1,1})\ne0$.
Suppose $(f_n^\pm)_*(\omega_{1,n-1,1})=0$ and hence $(f_n^\pm)_*(\omega_{n-1,1,1})=0$ as well. 
Then 
\[
(f_n^\pm)_*H^0(\Omega^2_{X_n/B})=
(f_n^\pm)_*F^2W_2H^2_\dR(U_n/B)=0
\]
by \eqref{ellK3-prop-2}.
Since $(f_n^\pm)_*$ is surjective onto 
$H^0(\Omega^2_{\cE_n^\pm/B})$, this 
contradicts with that $\cE^\pm_n$ are K3 surfaces (Proposition \ref{ellK3-prop1}).
This completes the proof.
\end{pf}

\medskip

We discuss the $p$-adic Beilinson conjecture for $K_3(\cE_n^-)$ with $n=3,4,6$.
Let $f_n^-:U_n\to\cE_n^-$ be the morphism \eqref{K3-1-eq5}
(recall that we only consider the cases $n=3,4,6$).
For $a\in \Q\setminus\{0,1\}$, we denote by $f_{n,a}:U_{n,a}\to\cE_{n,a}^-$
and $X_{n,a}$
the fibers at the closed point $t=a$ of $\Spec B$.
Let $\zeta_n$ be a primitive $n$-th root of unity.
Let $p>3$ be a prime at which $X_{n,a}$ has a good ordinary reduction.
Let $U_{a,n,\Q(\zeta_n)}=U_{a,n}\times_\Q\Spec\Q(\zeta_n)$ and
$U_{a,n,\Z_{(p)}[\zeta_n]}$ a smooth model over $\Z_{(p)}[\zeta_n]$.
Let
\begin{equation}
\Ross(\zeta_n)=\left\{\frac{1-x_0}{1-\zeta_nx_0},\frac{1-x_1}{1-\zeta_nx_1},\frac{1-x_2}{1+x_2}
\right\}\in K_3^M(\O(U_{a,n,\Z_{(p)}[\zeta_n]}))
\end{equation}
be the higher Ross symbol, which we think of being an element of Quillen's $K_3$.
Put
\[
\xi_n:=N_{\Q(\zeta_n)/\Q}(\Ross(\zeta_n))\in K_3(U_{n,a,\Z_{(p)}})
\]
where $N_{\Q(\zeta_n)/\Q}$ is the norm map in Quillen's $K$-theory.
Since $\Ross(\zeta_n)$ lies in the image of $K_3(X_{n,a,\Z_{(p)}[\zeta_n]})$ 
(\cite[Corollary 4.4]{As-Ross1}),
so does $\xi_n$, and hence there is a lifting $\wt{\xi_n}\in K_3(X_{n,a,\Z_{(p)}})$.
Put
\[
\omega_n:=(f_{n,a}^-)_*\omega_{1,n-1,1},\,
\eta_n:=(f_{n,a}^-)_*\eta_{1,n-1,1}\in H^2_\dR(\cE_{n,a}^-).
\]
Let
\[
\xymatrix{
\reg_\syn:K_3(\cE_{n,a,\Z_{(p)}}^-)\ar[r]&H^3_\syn(\cE_{n,a,\Z_{(p)}}^-,\Q_p(3))\cong
H^2_\dR(\cE_{n,a}^-/\Q)
}\]
be the syntomic regulator map where $\cE_{n,a,\Z_{(p)}}^-$ is a smooth integral model of 
$\cE_{n,a}^-$ over $\Z_{(p)}$.
Since $\eta_n$ is a unit root vector, 
one can show
that $\langle\reg_\syn(f^-_{n*}\wt{\xi_n}),\eta_n\rangle$ does not depend on the choice of lifting
in the same way as in \S \ref{ono-sect} (see the diagram \eqref{CD-Z}).
We have
\begin{align*}
\langle\reg_\syn((f_{n,a}^-)_*\wt{\xi_n}),\eta_n\rangle
&=\langle\reg_\syn(\wt{\xi_n}),(f^-_{n,a})^*\eta_n\rangle\\
&=\langle\reg_\syn(\Ross(\zeta_n)),(f^-_{n,a})^*\eta_n\rangle
+\langle\reg_\syn(\Ross(\zeta^{-1}_n)),(f^-_{n,a})^*\eta_n\rangle.
\end{align*}
Since
\[
\reg_\syn(\Ross(\zeta^{\pm1}_n))=
(1-\zeta_n)(1-\zeta^{-1}_n)\cF^{(\sigma)}_{\frac1n,\frac{n-1}{n},\frac12}(t)|_{t=a}(
\omega_{1,n-1,1}+\omega_{n-1,1,1})
+\text{(other terms)}
\]
by Theorem \ref{main-1}, we have
\begin{align*}
\langle\reg_\syn((f_{n,a}^-)_*\wt{\xi_n}),\eta_n\rangle
&=2(1-\zeta_n)(1-\zeta^{-1}_n)\cF^{(\sigma)}_{\frac1n,\frac{n-1}{n},\frac12}(t)|_{t=a}
\langle\omega_{1,n-1,1}+\omega_{n-1,1,1},(f^-_{n,a})^*\eta_n\rangle\\
&=2(1-\zeta_n)(1-\zeta^{-1}_n)\cF^{(\sigma)}_{\frac1n,\frac{n-1}{n},\frac12}(t)|_{t=a}
\times 2\langle\omega_n,\eta_n\rangle\quad \text{(by \eqref{ellK3-prop-3})}\\
&=4(1-\zeta_n)(1-\zeta^{-1}_n)\cF^{(\sigma)}_{\frac1n,\frac{n-1}{n},\frac12}(t)|_{t=a}
\langle\omega_n,\eta_n\rangle.
\end{align*}
Finally we note
\begin{align*}
\langle\omega_n,\eta_n\rangle
&=\langle(f^-_{n,a})^*(f^-_{n,a})_*\omega_{1,n-1,1},\eta_{1,n-1,1}\rangle\\
&=c\langle\omega_{1,n-1,1}+\omega_{n-1,1,1},\eta_{1,n-1,1}\rangle
\quad \text{(by \eqref{ellK3-prop-1})}\\
&=c\langle\omega_{n-1,1,1},\eta_{1,n-1,1}\rangle
\end{align*}
with $c\ne0$ a constant, and this does not vanish by Lemma \ref{unit-lem2} (3).
Summing up the above, we have the description of the $p$-adic regulator
for $K_3(\cE_{n,a}^-)$.
\begin{thm}\label{ellK3-thm}
Let $\sigma(t)=a^{1-p}t^p$ with $p>3$. 
Suppose that $\cE_{n,a}^-$ has a good ordinary reduction at $p$.
Then
\[
\frac{\langle\reg_\syn((f_{n,a}^-)_*\wt{\xi_n}),\eta_n\rangle}{\langle\omega_n,\eta_n\rangle}
=4(1-\zeta_n)(1-\zeta^{-1}_n)\cF^{(\sigma)}_{\frac1n,\frac{n-1}{n},\frac12}(t)|_{t=a}
\]
with $n=3,4,6$.
\end{thm}
\begin{conj}
Suppose that $\cE_{n,a}^-$ is a singular K3 surface over $\Q$.
Let $A_{n,a}=\sum a_nq^n$ be the corresponding Hecke eigenform of weight $3$, and
$\alpha_p$ the unit root of $T^2-a_pT+p^2$.
Then there is 
a constant $C_{n,a}\in \Q^\times$ not depending on $p$ such that
\[
L_p(A_{n,a},\omega_{\mathrm{Tei}}^{-1},0)=C_{n,a}(1-p^2\alpha^{-1}_p)\cF^{(\sigma)}_{\frac1n,\frac{n-1}{n},\frac12}(t)|_{t=a}.
\]
\end{conj}
Unlike the K3 surfaces in \cite{ono}, 
the author has not worked out on
the ($p$-adic) $L$-functions of $\cE^-_{n,a}$.

\end{document}